\documentclass[10pt,reqno]{amsart}
\usepackage{amsmath}
\usepackage{amsfonts}
\usepackage{amssymb}
\usepackage{graphicx}
\usepackage{amsthm,graphicx,color,yfonts}
\usepackage{pdfsync}
\usepackage{epstopdf}%
\usepackage{pifont}
\usepackage{bbm}
\usepackage{a4wide}

\usepackage[colorlinks=true]{hyperref}
\hypersetup{linkcolor=red,citecolor=blue,filecolor=dullmagenta,urlcolor=blue}

\usepackage{wrapfig}
\usepackage{tikz}
\usetikzlibrary{arrows,calc,decorations.pathreplacing}
\definecolor{light-gray1}{gray}{0.90}
\definecolor{light-gray2}{gray}{0.80}
\definecolor{light-gray3}{gray}{0.60}

\numberwithin{equation}{section}
\theoremstyle{plain}
\newtheorem{thm}{Theorem}[section]
\newtheorem{lem}[thm]{Lemma}
\newtheorem{prop}[thm]{Proposition}

\theoremstyle{definition}
\newtheorem{defn}{Definition}[section]

\theoremstyle{remark}
\newtheorem{rem}{Remark}[thm]

\newcommand{\R}{{\mathbb{R}}}
\newcommand{\Z}{{\mathbb{Z}}}

\newcommand{\La}{\widehat{L}_1}
\newcommand{\Lbk}{\widehat{L}^K_2}
\newcommand{\Lak}{\widehat{L}_1^K}
\newcommand{\Lb}{\widehat{L}_2}
\newcommand{\Ja}{\widehat{J}_1}
\newcommand{\Jbk}{\widehat{J}^K_2}
\newcommand{\Jak}{\widehat{J}_1^K}
\newcommand{\Jbb}{\widehat{J}^B_2}
\newcommand{\Jab}{\widehat{J}_1^B}
\newcommand{\Jb}{\widehat{J}_2}

\def\norm#1{\left\|#1\right\|}
\newcommand{\wh}{\widehat}
\newcommand{\vu}{\vec{u}}
\newcommand{\vv}{\vec{v}}
\newcommand{\vx}{\mathbf{x}}

\makeatother

\usepackage{babel}
\addto\shorthandsspanish{\spanishdeactivate{~<>}}

\addto\captionsenglish{}
\addto\captionsenglish{}
\addto\captionsenglish{}
\addto\captionsenglish{}
\addto\captionsenglish{}
\addto\captionsenglish{}
\addto\captionsenglish{}
\addto\captionsspanish{}
\addto\captionsspanish{}
\addto\captionsspanish{}
\addto\captionsspanish{}
\addto\captionsspanish{}
\addto\captionsspanish{}
\addto\captionsspanish{}

\begin{document}
	\title[Ill-posedness of abcd system]{Ill-posedness issues on $(abcd)$-Boussinesq system}

\author[C. Kwak]{Chulkwang Kwak}
\address{Department of Mathematics, Ewha Womans University, Seoul 03760, Korea}
\email{ckkwak@ewha.ac.kr}
\thanks{C. Kwak was partially supported by project France-Chile ECOS-Sud C18E06, and is partially supported by the Ewha Womans University Research Grant of 2020 and the National Research Foundation of Korea(NRF) grant funded by the Korea government(MSIT) (No. 2020R1F1A1A0106876811).
}

	\author[C. Maul\'en]{Christopher Maul\'en}  
	\address{Departamento de Ingenier\'{\i}a Matem\'atica and Centro
de Modelamiento Matem\'atico (UMI 2807 CNRS), Universidad de Chile, Casilla
170 Correo 3, Santiago, Chile.}
	\email{cmaulen@dim.uchile.cl}
	\thanks{Ch.M. was partially funded by Chilean research grants FONDECYT 1191412,  CONICYT PFCHA/DOCTORADO NACIONAL/2016-21160593 and CMM ANID PIA AFB170001.}

	\begin{abstract}
 In this paper, we consider the Cauchy problem for $(abcd)$-Boussinesq system posed on one- and two-dimensional Euclidean spaces. This model, initially introduced by Bona, Chen, and Saut \cite{BCS2002, BCS2004}, describes a small-amplitude waves on the surface of an inviscid fluid, and derived as a first order approximation of incompressible, irrotational  Euler equations. We mainly establish the ill-posedness of the system under various parameter regimes, which generalize the result of one-dimensional BBM-BBM case by Chen and Liu \cite{CL2013}. Most of results established here, we obtain the optimal result for two-dimensional BBM-BBM system.  The proof follows from an observation of the \emph{high to low frequency cascade} present in nonlinearity, motivated by Bejenaru and Tao \cite{BT2006}.
	\end{abstract}
	\maketitle

\section{Introduction}\label{sec:intro}
\subsection{Setting}

As a rigorous derivation from the free Eulerian formulation of water waves, Bona, Chen, and Saut \cite{BCS2004} proposed the model called one-dimensional $(abcd)$-Boussinesq, as
 \begin{equation}
 \label{eq:abcd}
 1D \; (abcd) \;\; \left\{ \begin{aligned}
 &(1-b\partial_{x}^{2})\partial_{t}\eta+\partial_{x}(a\partial_{x}^{2}u+u+u\eta) =0,\\
 &(1-d\partial_{x}^{2})\partial_{t}u+\partial_{x}(c\partial_{x}^{2}\eta+\eta+\frac{1}{2}u^{2}) =0,
 \end{aligned}\right.  \quad (t,x) \in \R \times \R.
 \end{equation}
As two-dimensional model, Bona, Colin and Lannes \cite{BCL2005}, formulated 2D $(abcd)$ as 
 \begin{equation}\label{eq:2D_ABCD}
2D \; (abcd) \;\; \left\{\begin{aligned}
&(1- b\,\Delta)\partial_t \eta + \nabla \cdot \left( a\, \Delta \vu + \vu + \vu \eta \right) =0,  \\
&(1- d\,\Delta)\partial_t \vu + \nabla \left( c\, \Delta \eta + \eta + \frac12 |\vu|^2 \right) =0,
\end{aligned}\right. \quad (t,\vx)\in \R\times\R^2.
\end{equation}
Here, unknowns $\eta$ and $u$ (also $\vu$) describe the free surface and the horizontal velocity of fluid, respectively. Both systems \eqref{eq:abcd} and \eqref{eq:2D_ABCD} are all first-order approximations of the incompressible and irrotational Euler equations assuming the small parameters defined by
\[
\alpha=\frac{A}{h} \ll 1, \quad
\beta=\frac{h^2}{\ell^2}\ll 1,\quad
\alpha\sim \beta,
\]
where $A$ and $\ell$ are typical wave amplitude and wavelength, and $h$ is the constant depth. Such assumptions sometimes referred to as small-amplitude long waves or Boussinesq or simply shallow water waves regimes (see \cite{Bous}). In the two-dimensional case, the irrotational hypothesis can be (mathematically) characterized as
\begin{equation}\label{eq:curl free}
\nabla \wedge \vu = 0,
\end{equation}
which is preserved by the evolution. Note that the condition \eqref{eq:curl free} is not necessary in the one-dimensional case since there is a single horizontal direction. See also \cite{AL2008} for relevant result.

\medskip

The parameters $(a,b,c,d)$ in both \eqref{eq:abcd} and \eqref{eq:2D_ABCD} are not arbitrary and hold the relations (see \cite{BCS2004})
\[\begin{gathered}
a=\frac12 \left(\theta^2-\frac13\right)\nu, \quad b=\frac12 \left(\theta^2-\frac13\right)(1-\nu),\\
c=\frac12\left(1-\theta^2\right)\nu -\tau ,\quad d=\frac{1}{2}(1-\theta^2)(1-\mu),
\end{gathered}\]
where $\theta \in [0,1]$ appears in the change of scaled horizontal velocity corresponding to the depth $(1-\theta)h$ below the undisturbed surface, $\tau$ is the surface tension ($\tau \geq0 $), and $\nu,\mu$ are arbitrary real numbers ensuring
\[a+b = \frac12\left(\theta^2 - \frac13\right), \quad c+d = \frac12(1-\theta^2) - \tau, \quad a+b+c+d = \frac13 - \tau.\]
The dispersive properties of the systems depend on the choice of the parameters. Precisely, the pair $(a,c)$ enhances the dispersion, while the pair $(b,d)$ weakens it (see \cite{CL2013}). This versatility makes the $(abcd)$-Boussinesq model interesting and challenging.

\medskip

Two systems 1D $(abcd)$ and $2D(abcd)$ allow the following energies
	\[E_{1D}[u,\eta](t)=\frac12 \int_{\R} (-au_x^2-c\eta_x^2+u^2(1+\eta)+\eta^2)(t,x)dx,\]
	and
		\begin{equation*} 
	E_{2D}[\vec{u},\eta](t)=\frac12 \int_{\R^2} (-a|\nabla \vec{u}|^2-c|\nabla \eta|^2+|\vec{u}|^2(1+\eta)+\eta^2)(t,x)dx,
	\end{equation*}
respectively, that both are conserved in time  when $b=d$ and $a,c <0$. Thus local well-posedness in $H^1$-level space is immediately extended to the global one at least for small data. Note that Sobolev embedding in two-dimensional case is not enough to control $L^{\infty}$ norm of $\eta$, but Gagliardo-Nirenberg interpolation inequality can control $\eta |\vec{u}|^2$.  

\medskip

These models have been extensively studied (in various perspective) in the literature, see e.g. \cite{BCS2002, BCS2004, B2006, DMS2007, LPS2011, SX2012, S_impa, CL2013, BCL2015, SWX2015, KMPP2019, MP2019, KM2020, SX2020,SX2020-1}. Among other them, we focus on Cauchy problems for these systems. In \cite{BCS2002, BCS2004}, Bona, Chen and Saut first studied local and global well-posedness of linear and nonlinear problems, and established the following results (the following results only exhibit the case when $\mathcal{H}$ (see \eqref{eq:h}) has order $0$):
\begin{enumerate}
\item the generic regime in $H^s(\R) \times H^s(\R)$, for $s\geq 0$.  
\item the BBM-BBM regime in $H^s(\R) \times H^s(\R)$, for $s\geq 0$. 
\item the KdV-KdV regime in $H^s(\R) \times H^s(\R)$, for $s> 3/4$. 
\end{enumerate}
In \cite{DMS2007}, Dougalis, Mitsotakis, and Saut proved that two-dimensional $(abcd)$ Boussinesq system under the generic regime is locally well-posed in $H^s(\R^2)\times H^s(\R^2)$ for $s>0$. Note that this local result is indeed valid in $L^2(\R^2) \times L^2(\R^2)$ by improving Grisvard's bilinear estimate \cite{Grisvard}, see Appendix \ref{app:A} (Lemma \ref{lem:bilinear}). 
In \cite{LPS2011}, Linares, Pilod, and Saut focused on the strongly dispersive (KdV-KdV system) regime, and established local well-posedness result in $H^s(\R^2)\times H^s(\R^2)$ for $s>3/2$. 
Previously, Schonbek \cite{S1981} and Amick \cite{A1984} considered a version of the original Boussinesq system ($a = c = b = 0, b =1/3$), and proved global well-posedness under a non-cavitation condition via parabolic regularization. Later, Burteau \cite{B2006} improved it without a non-cavitation condition. Studies on long time existence of solutions have been done in, for instance, \cite{SX2012,MSZ2012,SWX2015,SX2020,SX2020-1}. In these works, the authors established the well-posedness for large time with appropriate time scales.
 
\medskip
  
In contrast with results mentioned above, this paper concerns with the ill-posed issues on one- and two-dimensional $(abcd)$-Boussinesq systems in the following cases:
\begin{enumerate}
	\item Generic regime
\begin{equation}\label{eq:generic}
a,c<0, \quad b,d>0,
\end{equation}
\item KdV-KdV regime
\begin{equation}\label{eq:kdv_condition}
a=c=\frac16, \quad b=d=0,
\end{equation}
\item BBM-BBM regime
\begin{equation}\label{eq:bbm_condition}
a=c=0, \quad b=d=\frac16.
\end{equation}
\end{enumerate}
 
As far as the authors know, there is only few results for ill-posedness issues. Chen and Liu \cite{CL2013} established the (mild) ill-posedness result for one-dimensional system under the weakly dispersive regime (1D BBM-BBM system) below $L^2$.  The main idea follows the abstract theory developed by Bejenaru-Tao \cite{BT2006}. The authors also discussed the formation of singularities and provided blow-up criteria. Recently, \cite{ABM2019} Ambrose, Bona and Milgrom have established the ill-posedness of the one-dimensional periodic Kaup system ($a=1/3$ and $b=c=d=0$) is ill-posed in any positive regularity Sobolev space, in the sense that the flow map is discontinuous at the origin. They also concerned with the case that the generic condition \eqref{eq:generic} is negated.

\subsection{Main results}

Before presenting our results, we clarify what we mean ``ill-posedness". To do this, we first define ``well-posedness" of Partial Differential equations problems. As the author's best knowledge, The French mathematician Jacques Hadamard initially proposed the concept of well-posed problems as

\begin{defn}[Well-posedness]\label{def:well-posed}
The mathematical models of physical phenomena should have the following properties:

\bigskip

\begin{itemize}
\item there exists a solution,

\medskip

\item the solution is unique,

\medskip

\item the solution behaves continuously with the initial condition.
\end{itemize} 
\end{defn}

Obviously, problems that are not well-posed in the sense of Hadamard are termed ill-posed, in other words, invalidity of one of above properties makes problems be ill-posed. In this paper, in order to obtain ill-posedness results, we attack the third property in Definition \ref{def:well-posed}. A precise strategy follows the negation of Proposition \ref{thm:C^infty}.

\medskip

We are now ready to present our main theorems.
\begin{thm}\label{thm:1d}
		 The 1D-abcd system \eqref{eq:abcd} is ill-posed in the sense that the flow map from initial data to solutions  is discontinuous at the origin in $H^s(\R) \times H^s(\R)$, where 
\begin{enumerate}
	\item $s < -\frac12$ for the generic case (see \eqref{eq:generic}).\label{Thm1}
	\item $s < -\frac32$ for the KdV-KdV case (see \eqref{eq:kdv_condition}).\label{Thm2}
\end{enumerate}
\end{thm}
Analogously, 
\begin{thm}\label{thm:2d}
	The 2D-abcd system \eqref{eq:2D_ABCD} is ill-posed in the sense that the flow map from initial data to solutions  is discontinuous at the origin in $H^s(\R^2) \times H^s(\R^2)$, where
	\begin{enumerate}
	\item $s < -\frac12$ for the generic case (see \eqref{eq:generic}).\label{Thm4}
\item $s < -\frac32$ for the KdV-KdV case (see \eqref{eq:kdv_condition}) .\label{Thm5}
\item $s < 0$ for the BBM-BBM case (see \eqref{eq:bbm_condition}). \label{Thm3}
	\end{enumerate}
\end{thm}

\begin{rem}\label{rem:optimal}
The BBM-BMM case of the one-dimensional $(abcd)$-Boussinesq system has been dealt with by Chen and Liu \cite{CL2013}. However, the two-dimensional BBM-BBM system is considered here 	for the first time, and together with Appendix \ref{app:A},  we completely resolve Cauchy problem for it.
\end{rem}

The proofs of theorems follows the same idea developed by Bejenaru and Tao \cite{BT2006}, and motivated by an observation as follows: All nonlinear interactions are quadratic, thus \emph{high $\times$ high} interaction components over an appropriate short time depending on the frequency cause \emph{resonances} near the origin of  the resulting frequency. For this reason, the flow cannot disperse the high-frequency energy for this time so that the smoothness of the flow breaks below certain regularity. Note that this observation is simply applied to a one-dimensional problem, but it is non-trivial to construct initial data that can cause \emph{resonance} in two-dimensional case.

\medskip

It is easy to see that $(abcd)$ systems are completely coupled systems, thus an attempt at decoupling of (at least) the linear system must  take precedence in order to observe its propagators. Under generic regime, standard transforms (see \eqref{eq:transform1d} and \eqref{eq:CoV} for one- and two-dimensional cases, respectively) diagonalize the linear operator with eigenvalues $\sigma$ (see \eqref{eq:sigma}) for one-dimensional case and $\rho$ \eqref{eq:eigenvalues} for two-dimensional case of order $0$, while those under BBM-BBM and KdV-KdV regimes have order $-2$ and $2$, respectively. The difference of orders of eigenvalues directly affect the dispersive properties of  solutions, thus so flow maps of stronger dispersive systems can take rougher initial data. Such observations can be seen in, for instance, Lemma \ref{lem:cos sin}, and relevant lemmas.

\medskip

Our results are coherent with one-dimensional BBM-BBM system, generalized BBM equation and KdV equation. On one hand, In \cite{BBM_illposed}, the authors established that the flow map is not of class $C^2$, and warned that their result is not suitable to assert that the BBM-equation is ill-posed in $H^s$ for negative values of $s$. Conversely, in \cite{BBM_illposed_c}, the authors proved the discontinuity of  the flow map at the origin in $H^s$ for $s < 0$. On the other hand, our result for the KdV-KdV system differs from the ill-posedness result of the original KdV equation established by Molinet \cite{KdV_illposed}. The proof follows from an argument of functional analysis together with the discontinuity Miura transform and the validity of Kato smoothing effect of mKdV solutions. However, the same argument may not apply to the KdV-KdV system, since it has no such good structure. We also refer to, e.g., \cite{HM2013,K2009,KT2010,CP2019,P2011,CL2013} for relevant ill-posedness problems of single equations.

\medskip

In Appendix \ref{app:A}, we give a bilinear estimates of  $(1-\Delta)\nabla (fg)$, which slightly improve Grisvard's result \cite{Grisvard}. This improvement enable us to obtain the local well-posedness of two-dimensional $(abcd)$ system under generic and BBM-BBM regimes in $L^2(\R^2) \times L^2(\R^2) \times L^2(\R^2)$. As mentioned in Remark \ref{rem:optimal}, the well-posedness result for two-dimensional BBM-BBM system, in addition to Theorem \ref{thm:2d} \eqref{Thm3}, asserts the completion of Cauchy problem for it. The proof is based on Littlewood-Paley theory, and is completed by a delicate observation of frequency interactions.

\subsection*{Organization of the paper}	 
This paper is organized as follows: Section \ref{sec:2} devotes to introducing abstract and general well- and ill-posedness arguments developed by Bejenaru and Tao \cite{BT2006}, and to representing Boussinesq equations as linearly decoupled forms. In Sections \ref{sec:3} and \ref{sec:4} we prove Theorems \ref{thm:1d} and \ref{thm:2d}, respectively. In Appendix \ref{app:A} we briefly provide a refined bilinear estimate to establish the well-posedness of some classes of systems. In Appendices \ref{app:1D_decomposition}, \ref{app:C} and \ref{app:D}, we give precise computations for decomposition of quadratic terms in the nonlinearities.

\subsection*{Notations}	 
For $x,y \in \R_+$ ($=\R \cap (0, \infty)$), $x \lesssim y$ means that there exists $C>0$ such that $x \le Cy$, and $x \sim y$ means $x \lesssim y$ and $y \lesssim x$. For a Schwartz function $f$ in $x \in \R^d$, we denote  the Fourier transform of $f$ by $\mathcal F (f)$ or $\widehat{f}$ defined by
\[\widehat{f}(\xi)=\int _{\R^d} e^{-ix \cdot \xi}f(x) \;dx, \quad \xi \in \R^d,\]
and $\check{f}$ denotes the inverse Fourier transform of $f$ defined by
\[\check{f}(\xi)= \frac{1}{(2\pi)^d}\int _{\R^d} e^{ix \cdot \xi}f(\xi) \;d\xi, \quad x \in \R^d.\]
In the rest of sections, the following properties of the Fourier transform among others will be used frequently:
\[\mathcal F(f\ast g)(\xi) = \widehat{f}(\xi)\widehat{g}(\xi) \quad \mbox{and} \quad \mathcal F (\partial_{x_i}f) = i\xi_i \widehat{f}(\xi), \;\; 1\le i \le d.\]

\section{Preliminaries}\label{sec:2}
 
\subsection{Bejenaru-Tao's abstract theory}
Here, we briefly present the abstract well- and ill-posed theory initially introduced in \cite{BT2006}.  
 Consider the abstract equation
\begin{equation}\label{eq:abstract}
\vec{v} = \mathcal L(\vec{v}_0) + \mathcal N_{2}(\vec{v}, \vec{v}),
\end{equation}
where $\vec{v}_0 \in D$   is the initial data, and $\vec{v}\in S $ is a solution of abstract equation. Here,  $\mathcal L(\vec{v}_0)$  and $\mathcal N_{2}(\vec{v}, \vec{v})$ are the linear and bilinear part of Duhamel's formula, respectively.

In the context of this work, $\vec{v}$ is a solution to the $(abcd)$ system \eqref{eq:abcd} ($\vec{v} = (\eta, u)$) or \eqref{eq:2D_ABCD} ($\vec{v} = (\eta, \vu) = (\eta, u_1,u_2)$).
 First, we introduce the definition of quantitative well-posedness introduced in \cite{BT2006}.
\begin{defn}[Quantitative well-posedness, \cite{BT2006}]\label{def:well}
Let $(D,\|\cdot \|_{D})$ be a Banach space of initial data, and $(S,\|\cdot \|_{S})$ be a Banach space of space-time functions. We say that \eqref{eq:abstract} is quantitatively well-posed in $D,S$ if one has the estimate of the form
\[\|\mathcal L(\vec{v}_0)\|_{S} \lesssim \|\vec{v}_0\|_{D}\]
and
\begin{equation}\label{eq:nonlinearestimate}
 \|\mathcal N_{2}(\vec{v}, \vec{v})\|_{S} \lesssim \|\vec{v}\|_{S}^2,
 \end{equation}
for all $\vec{v}_0 \in D$ and $\vec{v} \in S$.
\end{defn}
As we can see in the introduction,  it is known that $abcd$-system (in one- or two-dimensional case) is locally well-posed in the sense of Definition \ref{def:well} (see \cite{BCS2004}).

\medskip

The following theorem asserts that the quantitative well-posedness indeed guarantees the analytic well-posedness, that is, the flow from given initial data to a solution is represented as a power series expansion of continuous nonlinear maps.

\begin{thm}[Theorem 3 in \cite{BT2006}]\label{thm:C^infty}
Suppose that \eqref{eq:abstract} is quantitatively well-posed in the space $D, S$. Then, there exist constants $C_0, \varepsilon_0 >0$ such that for all $\vec{v}_0 \in B_D(0,\varepsilon_0)$, there exists a unique solution $\vec{v}[\vec{v}_0] \in B_S(0, C_0\varepsilon_0)$ to \eqref{eq:abstract}. More specifically, if we define the non-linear maps $A_n : D \to S$ for $n \in \Z_{>0}$ by the recursive formula
\[\begin{aligned}
& A_1(\vec{v}_0) := \mathcal L(\vec{v}_0)\\
&A_n(\vec{v}_0) := \sum_{\substack{n_1, n_2 \ge 1 \\ n_1+n_2 =n}} \mathcal N_2(A_{n_1}(\vec{v}_0), A_{n_2}(\vec{v}_0)), \quad n >1,
\end{aligned}\]
then we have the absolutely convergent (in $S$) power series expansion
\[\vec{v}[\vec{v}_0] = \sum_{n=1}^{\infty} A_n(\vec{v}_0),\]
for all $\vec{v}_0 \in B_D (0,\varepsilon_0)$.
\end{thm}

On the other hand, Theorem \ref{thm:C^infty} alternatively says that one can prove ill-posedness of \eqref{eq:abstract}, once showing discontinuity of $A_n$, for some $n$, i.e., $A_n$ does not satisfy \eqref{eq:nonlinearestimate}. This observation can be precisely stated as follows: 

\begin{prop}[Proposition 1 in \cite{BT2006}]\label{prop:ill-posedness}
Suppose that \eqref{eq:abstract} is quantitatively well-posed in the Banach spaces $D$ and $S$, with a solution map $f \mapsto u[f]$ from a ball $B_D$ in $D$ to a ball $B_S$ in $S$. Suppose that these spaces are then given other norms $D'$ and $S'$, which are weaker than $D$ and $S$ in the sense that
\[\norm{\vec{v}_0}_{D'} \lesssim \norm{\vec{v}_0}_{D}, \quad \norm{\vec{v}}_{S'} \lesssim \norm{\vec{v}}_{S}.\]
Suppose that the solution map $\vec{v}_0 \mapsto \vec{v}[\vec{v}_0]$ is continuous from $(B_D, \norm{\;}_{D'})$ to $(B_S, \norm{\;}_{S'})$. Then for each $n$, the non-linear operator $A_n : D \to S$ is continuous from from $(B_D, \norm{\;}_{D'})$ to $(S, \norm{\;}_{S'})$.
\end{prop}

\subsection{Equivalent representation of $abcd$ systems}

This subsection devotes to rewriting $(abcd)$ systems \eqref{eq:abcd} and \eqref{eq:2D_ABCD} in the form of \eqref{eq:abstract}. We follow the arguments in \cite{BCS2004}  and \cite{DMS2007} for the one- and two-dimensional cases, respectively. 
\subsubsection{One-dimensional case.} We first deal with one-dimensional case. Applying the Fourier transform to the linear abcd system (\eqref{eq:abcd} without $u\eta$ and $\frac12 u^2$), we obtain
    \begin{equation}
        \frac{d}{dt}\begin{pmatrix}
\widehat{\eta} \\
\widehat{u}
\end{pmatrix}
+i\xi\begin{pmatrix}
0 & \omega_1(\xi)\\
\omega_2(\xi) & 0
\end{pmatrix}
\begin{pmatrix}
\widehat{\eta} \\
\widehat{u}
\end{pmatrix}=0, \label{eq:linear}
    \end{equation}
  where
  \begin{equation}\label{eq:omega}
  \omega_1(\xi):=\dfrac{1-a\xi^2}{1+b\xi^2} \quad \mbox{and} \quad\omega_2(\xi):=\dfrac{1-c\xi^2}{1+d\xi^2}.
  \end{equation}
Using the transform
\begin{equation}\label{eq:transform1d}\begin{pmatrix}
\eta\\
u
\end{pmatrix} = 
\begin{pmatrix}
\mathcal H & \mathcal H\\
1 & -1
\end{pmatrix}
\begin{pmatrix}
v \\
w
\end{pmatrix},\end{equation}
where $\mathcal{H}$ is the Fourier multiplier defined by
  \begin{equation}\label{eq:h}
    \widehat{\mathcal{H}f}(\xi)=h(\xi)\widehat{f}(\xi), \quad \mbox{with} \quad h(\xi):=\left(\dfrac{\omega_1 (\xi)}{\omega_2(\xi)}\right)^{\frac12}=\left(\dfrac{(1-a\xi^2)(1+d\xi^2)}{(1-c\xi^2)(1+b\xi^2)} \right)^{\frac12},
  \end{equation}
 the system \eqref{eq:linear} becomes a symmetric form as 
 \[
                \frac{d}{dt}\begin{pmatrix}
\widehat{v} \\
\widehat{w}
\end{pmatrix}
+i\xi\begin{pmatrix}
\sigma(\xi) & 0\\
0 & -\sigma(\xi)
\end{pmatrix}
\begin{pmatrix}
\widehat{v} \\
\widehat{w}
\end{pmatrix}=0, \]
  where
   \begin{equation}
   \sigma(\xi):=(\omega_1(\xi) \omega_2(\xi))^{\frac12}=\left( \dfrac{(1-a\xi^2)(1-c\xi^2)}{(1+b\xi^2)(1+d\xi^2)}\right)^{\frac12}.\label{eq:sigma}
   \end{equation}
 Note that \eqref{eq:omega} and \eqref{eq:h} are well-defined whenever the parameters $a,b,c$ and $d$ satisfy the generic or KdV-KdV cases.
Therefore, the system is written in the form
       \[
        \frac{\partial}{\partial t}\begin{pmatrix}
v \\
w
\end{pmatrix}+ B\begin{pmatrix}
v \\
w
\end{pmatrix}=0, \]
   where
   \[
   \begin{aligned}
   B\begin{pmatrix}
v \\
w
\end{pmatrix}
:=& \left(i\xi\begin{pmatrix}
\sigma(\xi) & 0\\
0 & -\sigma(\xi)
\end{pmatrix}\begin{pmatrix}
\widehat{v} \\
\widehat{w}
\end{pmatrix}\right)^{\vee},
    \end{aligned}
    \]
Coming back to the original variables $\eta$ and $u$, we write the linear system as
\[ \frac{\partial}{\partial t}\begin{pmatrix}
\eta \\
u
\end{pmatrix}+ A\begin{pmatrix}
\eta \\
u
\end{pmatrix}=0,\]
where an operator $A$ is determined by (also explicitly computed by taking the Fourier transform)
\[A = \begin{pmatrix}
\mathcal H & \mathcal H\\
1 & -1
\end{pmatrix}
B
\begin{pmatrix}
\mathcal H & \mathcal H\\
1 & -1
\end{pmatrix}^{-1}
,\] 

and thus the solutions to the linear $abcd$ system are of the form
    \[\begin{pmatrix}
\eta \\
u
\end{pmatrix}(t,x)=S(t)\begin{pmatrix}
\eta_0\\
u_0
\end{pmatrix}(x), \]
   where $S(t)$ is associated to the linear flow of the system generated by $A$. It is clear that $S(t)$ is a \textbf{unitary group} on $H^s(\R)\times H^s(\R)$ for any $s\in \R$. When
\[\omega_1(\xi)\omega_2(\xi) > 0,\]
the linear flow $S(t)$ can be expressed 
   \begin{equation}\label{eq:semigroup}
   \mathcal F \left(S(t) \begin{pmatrix} f \\ g\end{pmatrix}\right)
   :=\begin{pmatrix}
        \cos(\xi \sigma(\xi)t) &-i\sin(\xi \sigma(\xi)t)\dfrac{\omega_1(\xi)}{\sigma(\xi)}\\
        -i\sin(\xi \sigma(\xi)t)\dfrac{\omega_2(\xi)}{\sigma(\xi)} & \cos(\xi\sigma(\xi)t)
        \end{pmatrix}
        \begin{pmatrix}
        \widehat{f} \\ \widehat{g}
        \end{pmatrix}.
   \end{equation}
Note that
   \[
        \dfrac{\omega_1(\xi)}{\sigma(\xi)}
        = \left(\dfrac{\omega_1(\xi)}{\omega_2(\xi)}\right)^{\frac12} = h(\xi)
        \quad \mbox{and} \quad
        \dfrac{\omega_2
        (\xi)}{\sigma(\xi)}
        =   \left(\dfrac{\omega_2(\xi)}{\omega_1(\xi)}\right)^{\frac12} = \frac{1}{h(\xi)},
   \]
where $h$ is as in \eqref{eq:h}. Let
\begin{equation}\label{L}
\widehat{L}_1(t,\xi) = \cos(\xi \sigma(\xi)t) \quad \mbox{and} \quad \widehat{L}_2(t, \xi) = i\sin(\xi \sigma(\xi)t).
\end{equation}
Then we rewrite \eqref{eq:semigroup} as
   \begin{equation}\label{eq:semigroup_l}
\mathcal F  \left(S(t) \begin{pmatrix} f \\ g\end{pmatrix}\right)
:=\begin{pmatrix}
\widehat{L}_1(t,\xi)  &-h(\xi)\widehat{L}_2(t,\xi) \\
-(h(\xi))^{-1}\widehat{L}_2(t,\xi)  & \widehat{L}_1(t,\xi) 
\end{pmatrix}
\begin{pmatrix}
\widehat{f} \\ \widehat{g}
\end{pmatrix}.
\end{equation}
Duhamel's principle for the nonlinear system \eqref{eq:abcd} yields
\begin{equation}\label{eq:Duhamel}
\begin{aligned}
\begin{pmatrix} \eta \\ u \end{pmatrix} =&~{} S(t) \begin{pmatrix} \eta_0 \\ u_0 \end{pmatrix} - \int_0^t S(t-s) \partial_x\begin{pmatrix}  (1-b\partial_{x}^{2})^{-1}(\eta u) \\ (1-d\partial_{x}^{2})^{-1}(\frac12 u^2) \end{pmatrix} (s) \; ds\\
=:&~{} S(t) \begin{pmatrix} \eta_0 \\ u_0 \end{pmatrix} + \mathcal N_2\left(\begin{pmatrix} \eta \\ u\end{pmatrix}, \begin{pmatrix} \eta \\ u \end{pmatrix} \right).
\end{aligned}
\end{equation}

\bigskip

\noindent \textbf{KdV-KdV regime.}
Making a simple additional scaling, one may assume that  $a=c=1$, and obtain that the linear system of \eqref{eq:abcd}  satisfies
\[\frac{d}{dt}\begin{pmatrix}
\widehat{\eta} \\
\widehat{u}
\end{pmatrix}
+i\xi(1-\xi^2)\begin{pmatrix}
0 & 1\\
1 & 0
\end{pmatrix}
\begin{pmatrix}
\widehat{\eta} \\
\widehat{u}
\end{pmatrix}=0. \]
Analogously, we obtain

that the linear propagator $S_K(t)$ is represented as 
\begin{equation}\label{eq:semigroup_Kdv}
\mathcal F  \left(S_K(t) \begin{pmatrix} f \\ g\end{pmatrix}\right)
:=\begin{pmatrix}
\widehat{L^K_1}(t,\xi) &-\widehat{L^K_2}(t, \xi)\\
-\widehat{L^K_2}(t, \xi) & \widehat{L^K_1}(t,\xi)
\end{pmatrix}
\begin{pmatrix}
\widehat{f} \\ \widehat{g}
\end{pmatrix},
\end{equation}
where 
\begin{equation}\label{LK}
\widehat{L^K_1}(t,\xi) = \cos(\xi \sigma_K(\xi)t), \quad \widehat{L^K_2}(t, \xi) = i\sin(\xi \sigma_K(\xi)t)\quad  \mbox{and } \quad
\sigma_K(\xi)=1-\xi^2.
\end{equation}

\subsubsection{Two-dimensional case.}\label{ssec:2D}
We write the two-dimensional linear $abcd$ system (\eqref{eq:2D_ABCD} without $\vu \eta$ and $\frac12 |\vu|^2$) in the equivalent form (in the Fourier space, for fixed $\xi \in \R^2 \setminus \{\mathbf{0}\}$)
\begin{equation}\label{abcd_fourier}
\partial_t \begin{pmatrix} \wh{\eta} \\ \wh{u}_1 \\ \wh{u}_2 \end{pmatrix} +i|\xi| \mathcal A(\xi) \begin{pmatrix} \wh{\eta} \\ \wh{u}_1 \\ \wh{u}_2 \end{pmatrix}  =0,
\end{equation}
where
\begin{equation}\label{linearopmatrix}\mathcal A (\xi) = \begin{pmatrix}
0 & \frac{\xi_1}{|\xi|}\left(\frac{1-a|\xi|^2}{1+b|\xi|^2}\right) & \frac{\xi_2}{|\xi|}\left(\frac{1-a|\xi|^2}{1+b|\xi|^2}\right) \\
\frac{\xi_1}{|\xi|}\left(\frac{1-c|\xi|^2}{1+d|\xi|^2}\right) & 0 & 0 \\
\frac{\xi_2}{|\xi|}\left(\frac{1-c|\xi|^2}{1+d|\xi|^2}\right) & 0 & 0
\end{pmatrix}.\end{equation}
Define the Fourier symbol $\varrho(|\xi|)$ by
\begin{equation}\label{eq:eigenvalues}
\varrho(|\xi|) = \left(\frac{(1-a|\xi|^2)(1-c|\xi|^2)}{(1+b|\xi|^2)(1+d|\xi|^2)}\right)^{\frac12}.
\end{equation}
A straightforward computation yields that the matrix $\mathcal A(\xi)$ has three eigenvalues $0$, $\varrho(\xi)$ and $-\varrho(\xi)$, thus the matrix $\mathcal A(\xi)$ is diagonalized as follows:
\[\mathcal{P}^{-1}(\xi)\mathcal{A}(\xi)\mathcal{P}(\xi) = \begin{pmatrix}
0 & 0 & 0 \\
0 & \varrho(|\xi|) & 0 \\
0 & 0 & -\varrho(|\xi|)
\end{pmatrix},\]
where the \emph{block matrix} and its inverse are given by
\begin{equation}\label{eq:block matrices}
\mathcal{P}(\xi) = \begin{pmatrix}
0 & \varsigma(|\xi|) & -\varsigma(|\xi|) \\
-\frac{\xi_2}{|\xi|} & \frac{\xi_1}{|\xi|} & \frac{\xi_1}{|\xi|} \\
\frac{\xi_1}{|\xi|} & \frac{\xi_2}{|\xi|} & \frac{\xi_2}{|\xi|} \\
\end{pmatrix}
\end{equation}
and
\begin{equation}\label{eq:bolckinverse}
\mathcal{P}^{-1}(\xi) = \frac{1}{2\varsigma(|\xi|)}\begin{pmatrix}
0 & -2\varsigma(|\xi|)\frac{\xi_2}{|\xi|} & 2\varsigma(|\xi|)\frac{\xi_1}{|\xi|} \\
1 & \varsigma(|\xi|)\frac{\xi_1}{|\xi|} & \varsigma(|\xi|)\frac{\xi_2}{|\xi|} \\
-1 & \varsigma(|\xi|)\frac{\xi_1}{|\xi|} & \varsigma(|\xi|)\frac{\xi_2}{|\xi|} \\
\end{pmatrix},
\end{equation}
respectively, for
\begin{equation}\label{eq:varsigma}
	\varsigma(|\xi|) = \left(\frac{(1-a|\xi|^2)(1+d|\xi|^2)}{(1-c|\xi|^2)(1+b|\xi|^2)}\right)^{\frac12}.
	\end{equation}

We consider the following change of variables:
\begin{equation}\label{eq:CoV}
\begin{pmatrix} \wh{\mu} \\ \wh{\nu}_1 \\ \wh{\nu}_2 \end{pmatrix} = \mathcal{P}^{-1} (\xi) \begin{pmatrix} \wh{\eta} \\ \wh{u}_1 \\ \wh{u}_2 \end{pmatrix} = \begin{pmatrix} -\frac{\xi_2}{|\xi|}\wh{u}_1 + \frac{\xi_1}{|\xi|}\wh{u}_2 \\ \frac{\wh{\eta}}{2\varsigma(|\xi|)} +  \frac{\xi_1}{2|\xi|}\wh{u}_1 + \frac{\xi_2}{2|\xi|} \wh{u}_2 \\ -\frac{\wh{\eta}}{2\varsigma(|\xi|)} + \frac{\xi_1}{2|\xi|} \wh{u}_1 + \frac{\xi_2}{2|\xi|}\wh{u}_2 \end{pmatrix}.
\end{equation}
Note that $\widehat{\mu} = 0$, since the fluid is \emph{irrotational} (equivalently, \eqref{eq:curl free}). Thus, new variables $(\nu_1, \nu_2)$ finally determine an equivalent expression of the system \eqref{abcd_fourier}, and for \eqref{eq:2D_ABCD}, as
\[
\partial_t \begin{pmatrix}  \nu_1 \\ \nu_2 \end{pmatrix} + \mathcal B(-i\nabla) \begin{pmatrix} \nu_1 \\ \nu_2 \end{pmatrix}  =0, \quad \begin{pmatrix} \nu_1(t=0) \\ \nu_2(t=0) \end{pmatrix} = \begin{pmatrix} \nu_{1,0} \\ \nu_{2,0} \end{pmatrix},
\]
where $\mathcal B(-i\nabla)$ is the $2 \times 2$ matrix operator whose entries are pseudo-differential operators, with symbol
\[i|\xi|  \begin{pmatrix}
\varrho(|\xi|) & 0 \\
0 & -\varrho(|\xi|)
\end{pmatrix},\]
which is the skew-Hermitian matrix.

\medskip

\noindent Coming back to the original variables $\eta$, $u_1$ and $u_2$, we write the linear system as
\[ \frac{\partial}{\partial t}\begin{pmatrix}
\eta \\
u_1 \\
u_2
\end{pmatrix}+ \mathbf{A}\begin{pmatrix}
\eta \\
u_1 \\
u_2
\end{pmatrix}=0,\]
where the linear operator $\mathbf{A}$ is indeed given by the inverse Fourier transform of $i|\xi| \mathcal A(\xi)$ as in \eqref{abcd_fourier}. Thus, the solutions to the linear $abcd$ system are of the form
\[        \begin{pmatrix}
\eta \\
u_1\\
u_2
\end{pmatrix}(t,x)=\mathbf{S}(t)\begin{pmatrix}
\eta_0\\
u_{1,0}\\
u_{2,0}
\end{pmatrix}(x),\]
   where $\mathbf{S}(t)$ is associated to the linear flow of the system, generated by $\mathbf{A}$. When, $a,b,c$ and $d$ satisfy the generic condition, the linear flow $\mathbf{S}(t)$ can be expressed as, 
   \begin{equation}\label{eq:semigroup_2d}
\begin{aligned}
   \mathcal F  &\left(\mathbf{S}(t) \begin{pmatrix} f \\ g \\ h\end{pmatrix}\right)\\
   :=&~{}\mathcal{P}(\xi) \begin{pmatrix} 0 & 0 & 0\\ 0 & e^{i|\xi|\varrho(|\xi|)t} & 0\\0 & 0 & e^{-i|\xi|\varrho(|\xi|)t} \end{pmatrix} \mathcal{P}^{-1}(\xi)   \begin{pmatrix}
        \widehat{f} \\ \widehat{g} \\ \wh{h}
        \end{pmatrix}\\
        =&~{}  \begin{pmatrix} \widehat{J}_1(t,\xi) & \varsigma(|\xi|)\frac{i\xi_1}{|\xi|}\widehat{J}_2(t,\xi) & \varsigma(|\xi|)\frac{i\xi_2}{|\xi|}\widehat{J}_2(t,\xi) \\ 
        \frac{i\xi_1}{\varsigma(|\xi|)|\xi|}\widehat{J}_2(t,\xi) & \frac{\xi_1^2}{|\xi|^2}\widehat{J}_1(t,\xi) & \frac{\xi_1\xi_2}{|\xi|^2}\widehat{J}_1(t,\xi) 
        \\ \frac{i\xi_2}{\varsigma(|\xi|)|\xi|}\widehat{J}_2(t,\xi) & \frac{\xi_1\xi_2}{|\xi|^2}\widehat{J}_1(t,\xi) & \frac{\xi_2^2}{|\xi|^2}\widehat{J}_1(t,\xi) \end{pmatrix}   \begin{pmatrix}
        \widehat{f} \\ \widehat{g} \\ \wh{h}
        \end{pmatrix},
        \end{aligned}
   \end{equation}
   where
   \begin{equation}\label{eq:J2D}
   \begin{aligned}
   \widehat{J}_1(t,\xi)=\cos(|\xi|\varrho(|\xi|)t),\ \ \
   \widehat{J}_2(t,\xi)=\sin(|\xi|\varrho(|\xi|)t),
   \end{aligned}
   \end{equation}
   and $\varrho(|\xi|)$ is defined in \eqref{eq:eigenvalues}.

Duhamel's principle for the nonlinear system \eqref{eq:2D_ABCD} yields
\begin{equation}\label{eq:Duhamel_2d}
\begin{aligned}
\begin{pmatrix} \eta \\ u_1 \\ u_2 \end{pmatrix} =&~{} \mathbf{S}(t) \begin{pmatrix} \eta_0 \\ u_{1,0} \\ u_{2,0} \end{pmatrix} - \int_0^t \mathbf{S}(t-s) \begin{pmatrix} (1- b\,\Delta)^{-1}(\nabla \cdot (\eta \vu )) \\ \frac12 (1- d\,\Delta)^{-1}\partial_{x_1}|\vu|^2 \\ \frac12 (1- d\,\Delta)^{-1} \partial_{x_2}|\vu|^2 \end{pmatrix} (s) \; ds\\
=:&~{} \mathbf{S}(t) \begin{pmatrix} \eta_0 \\ u_{1,0} \\ u_{2,0} \end{pmatrix} + \mathcal N_2\left(\begin{pmatrix} \eta \\ u_1 \\ u_2\end{pmatrix}, \begin{pmatrix} \eta \\ u_1 \\ u_2 \end{pmatrix} \right).
\end{aligned}
\end{equation}

\begin{rem}
	Notice that the eigenvalues $\sigma$ and $\varrho$ (see  \eqref{eq:sigma} and \eqref{eq:eigenvalues}, respectively)  have the same radial behavior regardless of the dimension.
\end{rem}

\medskip

\noindent \textbf{Case $a=c$ and $b=d\geq0 $.}

In this case, we insert the conditions $a=c$ and $b=d$ into \eqref{linearopmatrix}, we then have 
\begin{equation}\label{abcd_fourier_ab}
\partial_t \begin{pmatrix} \wh{\eta} \\ \wh{u}_1 \\ \wh{u}_2 \end{pmatrix} +i|\xi| \varrho_{ab}(|\xi|) \mathcal A_{ab}(\xi) \begin{pmatrix} \wh{\eta} \\ \wh{u}_1 \\ \wh{u}_2 \end{pmatrix}  =0, \quad \xi \in \R^2 \setminus \{\mathbf{0}\}
\end{equation}
where
\begin{equation}\label{eq:eigenvalues_ab}
\mathcal A_{ab} (\xi) = \begin{pmatrix}
0 & \frac{\xi_1}{|\xi|} & \frac{\xi_2}{|\xi|} \\
\frac{\xi_1}{|\xi|} & 0 & 0 \\
\frac{\xi_2}{|\xi|}& 0 & 0
\end{pmatrix} 
\ \mbox{and }\ 
\varrho_{ab}(|\xi|) = \frac{1-a|\xi|^2}{1+b|\xi|^2}.
\end{equation}
Analogously as above (for the generic case), we can find three eigenvalues $0$, $1$ and $-1$ of the matrix $\mathcal A_{ab}(\xi)$. Thus the matrix $\mathcal A_{ab}(\xi)$ is diagonalized as follows:
\[\mathcal{P}_{ab}^{-1}(\xi)\mathcal{A}_{ab}(\xi)\mathcal{P}_{ab}(\xi) = \begin{pmatrix}
0 & 0 & 0 \\
0 & 1 & 0 \\
0 & 0 & -1
\end{pmatrix},\]
where the \emph{block matrix} $\mathcal{P}_{ab}$ and its inverse are represented as in \eqref{eq:block matrices} and \eqref{eq:bolckinverse}, respectively, with $\varsigma(|\xi|)=1$. Change the variables analogous to \eqref{eq:CoV}, then we have
\[\begin{pmatrix} \wh{\mu} \\ \wh{\nu}_1 \\ \wh{\nu}_2 \end{pmatrix} = \mathcal{P}_{ab}^{-1} (\xi) \begin{pmatrix} \wh{\eta} \\ \wh{u}_1 \\ \wh{u}_2 \end{pmatrix} = \begin{pmatrix} -\frac{\xi_2}{|\xi|}\wh{u}_1 + \frac{\xi_1}{|\xi|}\wh{u}_2 \\ \frac{\wh{\eta}}{2} +  \frac{\xi_1}{2|\xi|}\wh{u}_1 + \frac{\xi_2}{2|\xi|} \wh{u}_2 \\ -\frac{\wh{\eta}}{2} + \frac{\xi_1}{2|\xi|} \wh{u}_1 + \frac{\xi_2}{2|\xi|}\wh{u}_2 \end{pmatrix}.\]
As same as before,  $\widehat{\mu} = 0$, and new variables $(\nu_1, \nu_2)$  determine an equivalent expression of the system \eqref{abcd_fourier_ab} as
\[\partial_t \begin{pmatrix}  \nu_1 \\ \nu_2 \end{pmatrix} + \mathcal B_{ab}(-i\nabla) \begin{pmatrix} \nu_1 \\ \nu_2 \end{pmatrix}  =0, \quad \begin{pmatrix} \nu_1(t=0) \\ \nu_2(t=0) \end{pmatrix} = \begin{pmatrix} \nu_{1,0} \\ \nu_{2,0} \end{pmatrix},\]
where $\mathcal B_{ab}(-i\nabla)$ is the $2 \times 2$ matrix operator whose entries are pseudo-differential operators, with symbol
\[i|\xi| \varrho_{ab}(|\xi|) \begin{pmatrix}
1 & 0 \\
0 & -1
\end{pmatrix},\]
which is the skew-Hermitian matrix.

\medskip

Solutions to the linear $(abcd)$ system are of the form
\[\begin{pmatrix}
\eta \\
u_1\\
u_2
\end{pmatrix}(t,x)=\mathbf{S}_{ab}(t)\begin{pmatrix}
\eta_0\\
u_{1,0}\\
u_{2,0}
\end{pmatrix}(x),\]
where $\mathbf{S}_{ab}(t)$ is associated to the linear flow of the system, precisely expressed as, 
\begin{equation}\label{eq:semigroup_2d_ab}
\begin{aligned}
\mathcal F  &\left(\mathbf{S}_{ab}(t) \begin{pmatrix} f \\ g \\ h\end{pmatrix}\right)
=~{}  \begin{pmatrix} \widehat{J}^{ab}_1(t,\xi) & \frac{i\xi_1}{|\xi|}\widehat{J}^{ab}_2(t,\xi) & \frac{i\xi_2}{|\xi|}\widehat{J}^{ab}_2(t,\xi) \\ 
\frac{i\xi_1}{|\xi|}\widehat{J}^{ab}_2(t,\xi) & \frac{\xi_1^2}{|\xi|^2}\widehat{J}^{ab}_1(t,\xi) & \frac{\xi_1\xi_2}{|\xi|^2}\widehat{J}^{ab}_1(t,\xi) 
\\ \frac{i\xi_2}{|\xi|}\widehat{J}^{ab}_2(t,\xi) & \frac{\xi_1\xi_2}{|\xi|^2}\widehat{J}^{ab}_1(t,\xi) & \frac{\xi_2^2}{|\xi|^2}\widehat{J}^{ab}_1(t,\xi) \end{pmatrix}   \begin{pmatrix}
\widehat{f} \\ \widehat{g} \\ \wh{h}
\end{pmatrix},
\end{aligned}
\end{equation}
where
\begin{equation}\label{eq:J2D_ab}
\begin{aligned}
\widehat{J}^{ab}_1(t,\xi)=\cos(|\xi|\varrho_{ab}(|\xi|)t),\ \ \
\widehat{J}^{ab}_2(t,\xi)=\sin(|\xi|\varrho_{ab}(|\xi|)t),
\end{aligned}
\end{equation}
and $\varrho_{ab}(|\xi|)$ is defined in \eqref{eq:eigenvalues_ab}.

\medskip

Duhamel's principle for the nonlinear system \eqref{eq:2D_ABCD} yields
\begin{equation}\label{eq:Duhamel_2d_ab}
\begin{aligned}
\begin{pmatrix} \eta \\ u_1 \\ u_2 \end{pmatrix} =&~{} \mathbf{S}_{ab}(t) \begin{pmatrix} \eta_0 \\ u_{1,0} \\ u_{2,0} \end{pmatrix} - \int_0^t \mathbf{S}_{ab}(t-s) \begin{pmatrix} (1- b\,\Delta)^{-1}(\nabla \cdot (\eta \vu )) \\ 
\frac12 (1- b\,\Delta)^{-1}\partial_{x_1}|\vu|^2 \\ 
\frac12 (1- b\,\Delta)^{-1} \partial_{x_2}|\vu|^2 \end{pmatrix} (s) \; ds\\
=:&~{} \mathbf{S}_{ab}(t) \begin{pmatrix} \eta_0 \\ u_{1,0} \\ u_{2,0} \end{pmatrix} + \mathcal N_2\left(\begin{pmatrix} \eta \\ u_1 \\ u_2\end{pmatrix}, \begin{pmatrix} \eta \\ u_1 \\ u_2 \end{pmatrix} \right).
\end{aligned}
\end{equation}
\medskip

The above case covers KdV-KdV and BBM-BBM regimes, but we distinguish them below for simplicity.
\medskip

\noindent\textbf{KdV-KdV regime.}
As same as the one-dimensional case, one may assume that  $a=c=1$. Then, the system \eqref{eq:2D_ABCD}  has the form
\[\begin{cases}
&\partial_t \eta  + \nabla \cdot \left( \Delta \vu + \vu + \vu \eta \right) =0, \quad (t,\vx)\in \R\times\R^2, \\
&\partial_t \vu  + \nabla \left( \Delta \eta + \eta  + \frac12 |\vu|^2 \right) =0.
\end{cases}\]
Analogously, we have
\begin{equation}\label{eq:varrhoK}
\varrho_K(|\xi|) = 1-|\xi|^2.
\end{equation}
Due to $\varrho_K(|\xi|)\in \R$, the semigroup has the form
\[\mathcal{F} \left(\mathbf{S}_K(t)\right)=
\begin{pmatrix} \Jak(\xi,t) & \frac{i\xi_1}{|\xi|}\Jbk(\xi,t) & \frac{i\xi_2}{|\xi|}\Jbk(\xi,t) \\ \frac{i\xi_1}{|\xi|}\Jbk(\xi,t) & \frac{\xi_1^2}{|\xi|^2}\Jak(\xi,t) & \frac{\xi_1\xi_2}{|\xi|^2}\Jak(\xi,t) \\ \frac{i\xi_2}{|\xi|}\Jbk(\xi,t) & \frac{\xi_1\xi_2}{|\xi|^2}\Jak(\xi,t) & \frac{\xi_2^2}{|\xi|^2}\Jak(\xi,t) \end{pmatrix}.\]
where 
\begin{equation}\label{eq:J2D_K}
\begin{aligned}
\Jak(\xi,t)=\cos(|\xi|\varrho_K(|\xi|)t),\ \ \
\Jbk(\xi,t)=\sin(|\xi|\varrho_K(|\xi|)t).
\end{aligned}
\end{equation} 

\medskip

\noindent\textbf{BBM-BBM regime}
Let  $a=c=0$ and $b=d=1/6$. Then, the system \eqref{eq:2D_ABCD}  has the form 
\[\left\{
\begin{aligned}
\left(1- \frac16\,\Delta\right)\partial_t \eta  + \nabla \cdot \left( \vu + \vu \eta \right) = 0,  \\
\left(1- \frac16 \,\Delta\right)\partial_t \vu  + \nabla \left( \eta  + \frac12 |\vu|^2 \right) = 0,
\end{aligned}\right. \quad (t,\vx)\in \R\times\R^2.\]
Analogously, we have
\[\varrho_B(|\xi|) =\left(1+\frac16|\xi|^2\right)^{-1}. \]
Due to $\varrho_B\in \R$ and $\varsigma_B(|\xi|) =1$, the semigroup is 
\[\mathcal{F} \left(\mathbf{S}_B(t)\right)=
\begin{pmatrix} \Jab(\xi,t) & \frac{i\xi_1}{|\xi|}\Jbb(\xi,t) & \frac{i\xi_2}{|\xi|}\Jbb(\xi,t) \\
\frac{i\xi_1}{|\xi|} \Jbb(\xi,t) & \frac{\xi_1^2}{|\xi|^2}\Jab(\xi,t) & \frac{\xi_1\xi_2}{|\xi|^2}\Jab(\xi,t) \\
\frac{i\xi_2}{|\xi|}\Jbb(\xi,t) & \frac{\xi_1\xi_2}{|\xi|^2}\Jab(\xi,t) & \frac{\xi_2^2}{|\xi|^2}\Jab(\xi,t)
\end{pmatrix},\]
where
\begin{equation}\label{eq:J2D_BBM}
\begin{aligned}
\Jab(\xi,t)=\cos(|\xi|\varrho_B(|\xi|)t),\ \ \
\Jbb(\xi,t)=\sin(|\xi|\varrho_B(|\xi|)t).
\end{aligned}
\end{equation}

\section{Proof of Theorem \ref{thm:1d}}\label{sec:3}

\subsection{Generic regime: $a,~c < 0$ and $b,~d > 0$}

Before proving Theorem \ref{Thm1}, we address the  following two lemmas, which play key roles in our proof.
\begin{lem}\label{lem:sigma}
Let $\sigma(\xi)$ be as in \eqref{eq:sigma}. For $|\xi| \sim N$, we have the following relation:
\[\sigma(\xi)=\sqrt{\dfrac{ac}{bd}}+\tilde{\sigma}(\xi), \]
where $N$ is sufficiently large and $\tilde{\sigma}(\xi) = O(\xi^{-2})$ as $|\xi|\to \infty$.
\end{lem}
\begin{proof}
A straightforward computation gives
\[\begin{aligned}
 (\sigma(\xi))^2=~{}\dfrac{(1-a\xi^2)(1-c\xi^2)}{(1+b\xi^2)(1+d\xi^2)}
 =~{}\frac{ac}{bd} \left(1 + \frac{\alpha \xi^2 + \beta}{(1+b\xi^2)(1+d\xi^2)}\right),
\end{aligned}\]
where
\[\alpha = -(b+d) - \frac{bd(a+c)}{ac} \quad \mbox{and} \quad \beta = \frac{bd}{ac} - 1,\]
which implies
\[
\sigma(\xi)=\sqrt{\dfrac{ac}{bd}}\sqrt{1+ \dfrac{\alpha \xi^2+\beta}{(1+b\xi^2)(1+d\xi^2)}}.
\]
Using the binomial series expansion, we know
\[
(1+x)^{1/2}=\sum_{k=0}^{\infty} \binom{\frac12}{k}x^k=1+\frac{1}{2}x+O(x^2) \quad \mbox{for} \quad  |x|<1.
\]
 Thus, we can write
\[
\sigma(\xi)=\sqrt{\dfrac{ac}{bd}}+\tilde{\sigma}(\xi),
\]
where
\[
\tilde{\sigma}(\xi)= \sqrt{\dfrac{ac}{bd}} \cdot \dfrac{\alpha\xi^2 + \beta}{2(1+b\xi^2)(1+d\xi^2)} + O(\xi^{-4})\ \mbox{ as } |\xi| \gg 1,
\] 
since
\[\left| \dfrac{\alpha\xi^2+\beta}{(1+b\xi^2)(1+d\xi^2)} \right| \ll 1 \ \mbox{ for } |\xi| \gg 1.\]
This concludes the proof.
\end{proof}
Lemma \ref{lem:sigma} enables us to capture a specific nonlinear interaction among other interactions which makes non-smoothness of the flow, see Lemma \ref{lem:cos sin} below.
\begin{lem}\label{lem:cos sin}
Let $N \gg 1$ be sufficiently large, $T = \frac{1}{100N}$ and $0 \le s \le t \le T$. If $|\xi_1|, |\xi - \xi_1| \sim N$ and $|\xi| \sim 1$, then 
\[\left(\frac{1+d\xi^2}{1+b\xi^2}\right)\frac{h(\xi_1)}{h(\xi)}\wh{L}_2(t-s,\xi)\wh{L}_2(s,\xi_1)\wh{L}_1(s,\xi-\xi_1) + \dfrac{1}{2}\wh{L}_1(t-s,\xi)\wh{L}_1(s,\xi_1)\wh{L}_1(s,\xi-\xi_1) \ge \frac{1}{32},
\]
where $h$ and $\wh{L}_j$, $j=1,2$ are as in \eqref{eq:h} and \eqref{L}, respectively.
\end{lem}
\begin{proof}
A direct observation gives
\[\left|\frac{1+d\xi^2}{1+b\xi^2}\right| \le \max\left(1, \frac{d}{b}\right).\]
From the definition of $h$ and the sizes of $\xi$ and $\xi_1$, we immediately know
\[\left|\frac{h(\xi_1)}{h(\xi)}\right| \lesssim \max\left(1, \frac{ac}{bd}\right).\]
Moreover, we also know from Lemma \ref{lem:sigma} that
\[\xi \sigma(\xi) = \sqrt{\dfrac{ac}{bd}}\xi + O(\xi^{-1}),\]
as $|\xi| \to \infty$. On one hand, the conditions $|\xi| \sim 1$ and $0 \le s \le t \le T$ with $T = \frac{1}{100N}$ yield
\[|\sin(\xi\sigma(\xi)(t-s))| \lesssim \frac1N,\] 
which implies
\[\frac{h(\xi_1)}{h(\xi)}\left|\wh{L}_2(t-s,\xi)\wh{L}_2(s,\xi_1)\wh{L}_1(s,\xi-\xi_1)\right|  \le \frac{1}{32},\]
for sufficiently large $N$. On the other hand, since
\[|\xi|\sigma(\xi)|t-s|, |\xi_1|\sigma(\xi_1)|s|, |\xi-\xi_1|\sigma(\xi-\xi_1)|s| \le \frac{\pi}{3},\]
for sufficiently large $N$, we obtain
\[\dfrac{1}{2}\wh{L}_1(t-s,\xi)\wh{L}_1(s,\xi_1)\wh{L}_1(s,\xi-\xi_1) \ge \frac{1}{16}.\]
Collecting all, we complete the proof of Lemma \ref{lem:cos sin}.
\end{proof}

\begin{proof}[\textbf{Proof of Theorem} \ref{thm:1d} \eqref{Thm1}]

We use a contradiction argument. Suppose that the flow map $\vv_0 \mapsto \vv[\vv_0]$ is continuous in $H^s$, $s < -\frac12$. Then, from Proposition \ref{prop:ill-posedness}, the map $\vv_0 \mapsto A_2(\vv_0)$ is also continuous, where
\[A_2(\vec{v}_0)=\mathcal N_2 (A_1(\vv_0),A_1(\vv_0)).\]
In what follows, we are going to prove that the map $\vv_0 \mapsto A_2(\vv_0)$ violates the following inequality:
\[
\|A_2(\vv_0)\|_{X^{s'}_{T}}\lesssim \| \vv_0\|_{s}^2,
\]
for $s < -\frac12$ and $s' \in \R$, which completes the proof.

\medskip

Let $\vv_0$ be an initial datum, which will be chosen later. Using \eqref{eq:semigroup_l}, we write
\begin{equation}\label{eq:v1}
\vv_{1}=\begin{pmatrix} \eta_1 \\ u_1  \end{pmatrix} =S(t) \begin{pmatrix} \eta_0 \\ u_0  \end{pmatrix}
 =\begin{pmatrix} L_1 \eta_0-h(-i\partial_x) L_2  u_0  \\ -(h(-i\partial_x))^{-1} L_2  \eta_0 +L_1 u_0  \end{pmatrix}
,\end{equation}
where $h$ and $L_j$, $j=1,2$, are as in \eqref{eq:h} and \eqref{L}, respectively. Let
\begin{equation}\label{eq:A_2}
A_2(\vec{v_0})=\int_{0}^t S(t-s)\partial_x \begin{pmatrix}
(1-b\partial_{x}^2)^{-1}(\eta_1u_1)\\
(1-d\partial_{x}^2)^{-1} \dfrac{u_1^2}{2}
\end{pmatrix}(s) \; ds =: \int_{0}^t \begin{pmatrix}
    Q_1(s) \\
     Q_2(s)
\end{pmatrix} \; ds,
\end{equation}
as in \eqref{eq:Duhamel}.

\medskip

For $N$ large enough, we choose the initial data $\eta_0$ as the zero function and $u_0$ as a large frequency localized function, more precisely, $\vec{v}_{0}=(\eta_0,u_0)$ so that
\begin{equation}\label{eq:initial_data}
    \widehat{\eta}_{0}=0 \quad \mbox{and} \quad \widehat{u}_{0}(\xi)=N^{-s}\chi_{\mathcal A_N}(\xi),
\end{equation}
where $\chi$ is the characteristic function and the set $\mathcal A_N$ is given by
\[\mathcal A_N = \left\{\xi \in \R : N-\frac{1}{2}\leq |\xi|\leq N+\frac{1}{2}\right\}.\]
Note that $\norm{u_0}_{H^s} \sim 1$. 
Inserting the initial data \eqref{eq:initial_data} into \eqref{eq:v1}, thus \eqref{eq:A_2}, we obtain $Q_2 = Q_{21} + Q_{22}$, where
\[\widehat{Q}_{21} =  \frac{i \xi}{h(\xi)(1+b\xi^2)}\int_{\R} h(\xi_1) \Lb(t-s,\xi) \Lb(s,\xi_1) \La(s,\xi-\xi_1) \widehat{u_0}(\xi_1) \widehat{u_0}(\xi-\xi_1)d\xi_1\]
and
\[\widehat{Q}_{22} = \frac{i \xi}{2(1+d\xi^2)} \int_{\R} \La(t-s,\xi) \La(s,\xi_1) \La(s,\xi-\xi_1)
\widehat{u_0}(\xi_1) \widehat{u_0}(\xi-\xi_1)d\xi_1.\]
A precise computation is given in Appendix \ref{app:1D_decomposition}. From the supports of $\widehat{u}_{0}(\xi_1)$ and $\widehat{u}_{0}(\xi-\xi_1)$, the possible values of $\xi$ satisfy
\[2N-1\leq |\xi|\leq 2N+1 \quad \mbox{or} \quad |\xi|\leq 1.\]
Moreover, by Lemma \ref{lem:sigma}, we have
\[
|\xi_1\sigma(\xi_1)|,|(\xi-\xi_1)\sigma(\xi-\xi_1)|\sim \sqrt{\dfrac{ac}{bd}}N+O(N^{-1}),
\]
for sufficiently large $N$.  Set $t:= \frac{1}{100N}$. Then, from Lemma \ref{lem:cos sin}, we obtain
\[ \begin{aligned}
        \Vert A_2 (\vec{v}_0) \Vert_{H^{s'}(\R)\times H^{s'}(\R)}
        &\ge \left\Vert \langle \xi\rangle^{s'}\int_0^t  \begin{pmatrix} \widehat{Q}_1 \\ \widehat{Q}_2 \end{pmatrix} \; ds \right\Vert_{(L^2 \times L^2)(|\xi|\le 1)}\\
        &\ge \left\Vert \langle \xi\rangle^{s'} \int_0^t  \widehat{Q}_2 ds \right\Vert_{L^2(|\xi| \le 1)}\\
         &\geq \frac{1}{32}\left\Vert \langle \xi\rangle^{s'} \dfrac{i\xi}{1+d\xi^2} \int_0^t  \int_{\R}\widehat{u}_{0}(\xi_1)\widehat{u}_{0}(\xi-\xi_1) \; d\xi_1  ds \right\Vert_{L^2(|\xi|\le 1)}\\
        &\gtrsim N^{-2s-1}\|\xi\|_{L^2(|\xi|\le 1)},
    \end{aligned}\]
which does not guarantee the uniform boundedness of $\Vert A_2 (\vec{v}_0) \Vert_{H^{s'}(\R)\times H^{s'}(\R)}$ for $s < -\frac12$.  This complete the proof.
\end{proof}

\subsection{KdV-KdV regime: $b= d = 0$ and $a=c=\frac16$}
In order to prove Theorem \ref{thm:1d} \eqref{Thm2}, we need a modification of Lemma \ref{lem:cos sin}, which was essential to prove the generic case. We have
\begin{lem}\label{lem:cos sin_KdV}
	Let $N \gg 1$ be sufficiently large and $T = \frac{1}{100N^3}$ and $0 \le s \le t \le T$. If $|\xi_1|, |\xi - \xi_1| \sim N$ and $|\xi| \sim 1$, then we have
	\[
	\begin{aligned}
			\wh{L^K_2}(t-s,\xi)\wh{L^K_2}(s,\xi_1)\wh{L^K_1}(s,\xi-\xi_1)
			+\dfrac{1}{2}\wh{L^K_1}(t-s,\xi)\wh{L^K_1}(s,\xi_1)  \wh{L^K_1}(s,\xi  & -\xi_1) \ge \frac{1}{32},
	\end{aligned}
	\]
	where $\sigma_K$ and $\wh{L^K_j}$, $j=1,2$ are as in \eqref{LK}, respectively.
\end{lem}

\begin{proof}
The proof is analogous to the proof of Lemma \ref{lem:cos sin}, thus we omit the details.
\end{proof}

\begin{proof}[\textbf{Proof of Theorem} \ref{thm:1d} \eqref{Thm2}]
	The proof follows the proof of Theorem \ref{thm:1d} \eqref{Thm1}. Suppose that the flow map $\vv_0 \mapsto \vv[\vv_0]$ is continuous in $H^s$, $s < -3/2$.

	\medskip

	Let $\vv_0=(\eta_0,u_0)$ be initial data given by \eqref{eq:initial_data}, and
	\[\vv_{1}=\begin{pmatrix} \eta_1 \\ u_1  \end{pmatrix} =S_K(t) \begin{pmatrix} \eta_0 \\ u_0  \end{pmatrix},\]
	where $S_K(t)$ is given by \eqref{eq:semigroup_Kdv}. Then, we write
	\[
	A_2(\vec{v_0})=\int_{0}^t S_K(t-s)\partial_x \begin{pmatrix}
	\eta_1u_1\\
	 \dfrac{u_1^2}{2}
	\end{pmatrix}(s) \; ds =: \int_{0}^t \begin{pmatrix}
	Q_1(s) \\
	Q_2(s)
	\end{pmatrix} \; ds.
	\]
	A straightforward computation enables us to decompose $Q_1$ and $Q_2$ as $Q_1 = Q_{11} + Q_{12}$ and $Q_2 = Q_{21} + Q_{22}$, where
	\[\widehat{Q}_{11} =-i \xi \int_{\R} \Lak(t-s,\xi) \Lbk(s,\xi_1)\Lak(s,\xi-\xi_1) \widehat{u_0}(\xi_1) \widehat{u_0}(\xi-\xi_1)d\xi_1,\]
	\[\widehat{Q}_{12} = - \frac{i \xi}{2} \int_{\R} \Lbk(t-s,\xi) \Lak(s,\xi_1)  \Lak(s,\xi-\xi_1) \widehat{u_0}(\xi_1) \widehat{u_0}(\xi-\xi_1)d\xi_1,\]
	\[\widehat{Q}_{21} = i \xi\int_{\R}  \Lbk(t-s,\xi)\Lbk(s,\xi_1) \Lak(s,\xi-\xi_1) \widehat{u_0}(\xi_1) \widehat{u_0}(\xi-\xi_1)d\xi_1\]
	and
	\[\widehat{Q}_{22} = \frac{i \xi}{2} \int_{\R}\Lak(t-s,\xi) \Lak(s,\xi_1) \Lak(s,\xi-\xi_1) \widehat{u_0}(\xi_1) \widehat{u_0}(\xi-\xi_1)d\xi_1,\]
	where $L_j^K$, $j=1,2$, are given in \eqref{LK}.
	On the supports of $\widehat{u}_{0}(\xi_1)$ and $\widehat{u}_{0}(\xi-\xi_1)$, the resulting frequency $\xi$ possibly lies in the regions
	\[2N-1\leq |\xi|\leq 2N+1 \quad \mbox{or} \quad |\xi|\leq 1.\]
	Moreover, for sufficiently large $N$, we observe
	\[
	|\xi_1\sigma_K(\xi_1)|,|(\xi-\xi_1)\sigma_K(\xi-\xi_1)|\sim N^3.
	\]
	Set $t:= \frac{1}{100N^3}$. From Lemma \ref{lem:cos sin_KdV}, we conclude
	\[ \begin{aligned}
	\Vert A_2 (\vec{v}_0) \Vert_{H^{s'}(\R)\times H^{s'}(\R)}
	&\ge  \left\Vert \int_0^t  \widehat{Q}_2 ds \right\Vert_{L^2(|\xi| \le 1)}\\
	&\geq \frac{1}{32}\left\Vert \langle \xi\rangle^{s'} i\xi \int_0^t  \int_{\R}\widehat{u}_{0}(\xi_1)\widehat{u}_{0}(\xi-\xi_1)  \;\; d\xi_1ds \right\Vert_{L^2_{\xi}(|\xi|\le 1))}\\
	&\gtrsim N^{-2s-3}\|\xi\|_{L^2(|\xi|\le 1)},
	\end{aligned}\]
	which does not guarantee the uniform boundedness of $\Vert A_2 (\vec{v}_0) \Vert_{H^{s'}(\R)\times H^{s'}(\R)}$ for $s < -\frac32$ and this completes the proof.
\end{proof}

\section{Proof of Theorem \ref{thm:2d}}\label{sec:4}

\subsection{Generic regime: $a,~c<0$ and $b,~d>0$}
We first address the  following lemma, which plays a similar role as Lemma \ref{lem:cos sin}.
\begin{lem}\label{lem:cos sin_2D}
	Let $N \gg 1$ be sufficiently large, and $T = \frac{1}{100N}$ and $0 \le s \le t \le T$. If $|\kappa|, |\xi - \kappa| \sim N$, $|\xi| \sim 1$ and $(\kappa-\xi)\cdot \kappa >0$ , then we have
	\begin{equation*}
	\begin{aligned}
	\frac{(1+d|\xi|^2)}{(1+b|\xi|^2)} 
	\frac{\xi\cdot \kappa}{|\xi| |\kappa| }
	&\frac{\varsigma(|\xi-\kappa|)}{ \varsigma(|\xi|) } \Jb(t,\xi)  \Jb(s,\xi-\kappa) \Ja(s,\kappa)\\
	&+
	\frac{(\kappa-\xi)\cdot \kappa}{|\xi-\kappa| |\kappa|}
	\frac12 \Ja(t, \xi) \Ja(s,\xi-\kappa)\Ja(s,\kappa) 
	\geq \frac{1}{16}\left( 	\frac{(\kappa-\xi)\cdot \kappa}{|\xi-\kappa||\kappa|} -\frac{1}{2}\right)
	\end{aligned}
	\end{equation*}
	where $\cdot$ denotes the standard inner product in Euclidean space, and $\varsigma$, $\wh{J}_j$, $j=1,2$  are as in \eqref{eq:eigenvalues}, \eqref{eq:J2D}, respectively.
\end{lem}

\begin{proof}
	A direct computation yields
	\begin{equation*} 
	\left|\dfrac{\varsigma(|\xi-\kappa|)}{\varsigma(|\xi|)} \right|
	\leq \max\left( 1,\dfrac{ab}{cb}\right)
	\ \mbox{and} \ \ \
	\left|\frac{\xi\cdot \kappa}{|\xi| |\kappa| } 
	\frac{(1+d|\xi|^2) }{(1+b|\xi|^2) }
	\right|
	\leq \max\left(1,\frac{d}{b}\right).
	\end{equation*}
	Note that Lemma \ref{lem:sigma} is valid for $\varrho$, thus
	\[|\kappa| \varrho(|\kappa|) = \sqrt{\dfrac{ac}{bd}}|\kappa| + O(|\kappa|^{-1}),\]
	as $|\kappa| \to \infty$. On one hand, since $|\xi| \sim 1 $ and $0 \le s \le t \le T$ with $T = \frac{1}{100N}$, we know $|\sin(|\xi|\varrho(|\xi|)(t-s))| \lesssim \frac{1}{N}$, hence
	
	\[
	\begin{aligned}
	\bigg|
	\frac{\xi\cdot \kappa}{|\xi| |\kappa| } \frac{\varsigma(|\xi-\kappa|)}{\varsigma(|\xi|)}
	\frac{(1+d|\xi|^2) }{(1+b|\xi|^2) }
	\Jb(t-s,\xi) \Jb(s,\xi-\kappa) \Ja(s,\kappa) 
	\bigg|\leq \frac{1}{32}
	\end{aligned}
	,\]
	for $N$ large enough. On the other hand, since
	\[|\xi|\varrho(|\xi|)|t-s|, |\kappa|\varrho(|\kappa|)|s|, |\xi-\kappa|\varrho(|\xi-\kappa|)|s| \le \frac{\pi}{3},\]
	we have
	\[
	\frac{1}{2}\wh{J}_1(t-s,\xi)\wh{J}_1(s,\kappa)\wh{J}_1(s,\xi-\kappa) \ge \frac{1}{16}.
	\]
	We complete the proof from the last frequency condition
		\[
     \frac{(\kappa-\xi)\cdot \kappa }{|\kappa-\xi||\xi|}>0.
	\]
\end{proof}
\begin{rem}
The last condition $(\kappa-\xi)\cdot \kappa >0$ in Lemma \ref{lem:cos sin_2D} is not artificial under the rest conditions $|\kappa|, |\xi - \kappa| \sim N$ and $|\xi| \sim 1$, since the low resulting frequency from two high frequencies interaction occurs only when two high frequencies lie in the opposite side around the origin. A precise observation will be seen in the proof of Theorem \ref{thm:2d} below.
\end{rem}
\begin{proof}[\textbf{Proof of Theorem} \ref{thm:2d} \eqref{Thm4}]

Analogously to Theorem \ref{thm:1d} \eqref{Thm1}, suppose that the flow map $\vv_0 \mapsto \vv[\vv_0]$ is continuous in $H^s$, $s < -\frac12$.
\medskip

Let $\vv_0=(\eta_0,u_{01},u_{02})$ be an initial data to be chosen later. Let
\begin{equation}\label{eq:vv_1}
\vv_{1}=\begin{pmatrix} \eta_1 \\ u_{11}\\u_{12}  \end{pmatrix} =\mathbf{S}(t) \begin{pmatrix} \eta_0 \\ u_{01}\\ u_{02}  \end{pmatrix}
=
\left(\begin{array}{c}
J_1\eta_0+\varsigma(|\nabla|)\frac{\partial_{x_1}}{|\nabla|} J_2 u_{01}+\varsigma(|\nabla|)\frac{\partial_{x_2}}{|\nabla|} J_2 u_{02} \\
\frac{\partial_{x_1}}{\varsigma(|\nabla|)|\nabla|} J_2  \eta_0+\frac{-\partial_{x_1}^2}{|\nabla|^2}J_1 u_{01} +\frac{-\partial{x_1}\partial_{x_2}}{|\nabla|^2}J_1 u_{02}\\
\frac{\partial_{x_2}}{\varsigma(|\nabla|)|\nabla|} J_2  \eta_0 +\frac{-\partial_{x_1}\partial_{x_2}}{|\nabla|^2}J_1 u_{01}+\frac{-\partial_{x_2}^2}{|\nabla|^2}J_1 u_{02}
\end{array}\right),
\end{equation}
where $J_1, J_2$ are defined in \eqref{eq:J2D}, and $\varsigma(|\xi|)$ is given in \eqref{eq:varsigma}.

\medskip

For $N$ large enough, set

\[\begin{aligned}
\mathcal{S}_N := &\left\{ \kappa \in \R^2 :N-\frac12 \leq \kappa_1 \leq N+\frac12 \mbox{ and } |\kappa_2|\leq 1 \right\}\\
&\cup \left\{ \kappa \in \R^2 :-N-\frac12 \leq \kappa_1\leq -N+\frac12 \mbox{ and } |\kappa_2|\leq 1 \right\}
=:\mathcal{S}_N^{+}\cup \mathcal{S}_N^{-}.
\end{aligned}\]
We choose the initial data $\vec{v}_{0}=(\eta_0,u_{01},u_{02})$ as
\begin{equation}\label{eq:initial_data_2D}
\eta_0=0, \quad \mbox{and} \quad u_{01} = u_{02}=\psi_N,
\end{equation}
where $\widehat{\psi}_N(\xi)=N^{-s}\chi_{\mathcal S_N}(\xi)$. Note that $\| \psi_N\|_{H^s} \sim 1$. Moreover, on the supports of $\widehat{\psi}_N(\xi-\kappa)$ and $\wh{\psi}_N(\kappa)$, the resulting frequency $\xi$ belongs to
\[S_L:=\{ \xi\in \R^2 : |\xi_j|\leq 2, j=1,2\}\]
or
\[S_H:= \left\{ \xi\in \R^2 :
2N-1 \leq |\xi_1|\leq 2N+1,\ \ 
|\xi_2|\leq 2\right\} .\]
From \eqref{eq:Duhamel_2d}, we have
\begin{equation}\label{eq:Duhamel_2d_1}
A_2(\vec{v_1})=\int_{0}^t \mathbf{S}(t-s)
    \left(\begin{array}{c}
                (1-b\Delta)^{-1}[\partial_{x_1}(\eta_{1} u_{11})+\partial_{x_2}(\eta_1 u_{12})]     \\
                2^{-1}(1-d\Delta)^{-1}\partial_{x_1}(u_{11}^2+ u_{12}^2) \\
                2^{-1}(1-d\Delta)^{-1}\partial_{x_2}(u_{11}^2 +u_{12}^2)
    \end{array}\right)
    =: \int_{0}^t \left(\begin{array}{c}
     Q_1(s) \\
     Q_2(s)\\
     Q_3(s)
    \end{array} \right) ds.
    \end{equation}
    Inserting \eqref{eq:initial_data_2D} into \eqref{eq:vv_1}, thus \eqref{eq:Duhamel_2d_1}, we obtain $\widehat{Q}_2$ as $Q_2 = Q_{21} + Q_{22}$, where 
\[\begin{aligned}
\wh{Q}_{21}=&-
i\xi_1 \int_{\R^2} \frac{p(\xi,\kappa)}{(1+b|\xi|^2)} 
\frac{\xi\cdot \kappa}{|\xi| |\kappa| }
\frac{\varsigma(|\xi-\kappa|)}{ \varsigma(|\xi|) }  \Jb(t,\xi)  \Jb(s,\xi-\kappa) \Ja(s,\kappa) \widehat{\psi}_N(\xi-\kappa) \widehat{\psi}_N(\kappa) d\kappa,\\
\wh{Q}_{22}=& -i \xi_1  
\int_{\R^2} \frac{(\kappa-\xi)\cdot \kappa}{|\xi-\kappa| |\kappa|}
\frac{p(\xi,\kappa)}{(1+d|\xi|^2)}\frac12 \Ja(t, \xi) \Ja(s,\xi-\kappa)\Ja(s,\kappa) \widehat{\psi}_N(\xi-\kappa)  \widehat{\psi}_N (\kappa)d\kappa,
\end{aligned}
\]
for
\begin{equation}\label{eq:pp}
p(\xi,\kappa)=\frac{(\xi_1+\xi_2-\kappa_1-\kappa_2)}{|\xi-\kappa|}
\frac{(\kappa_1+\kappa_2)}{|\kappa|}.
\end{equation}
See Appendix \ref{app:C} for precise computations of $Q_j$, $j=1,2,3$.  

For $\xi-\kappa,~ \kappa \in \subset \mathcal{S}_N$, the resulting frequency $\xi$ lies in $\mathcal S_L$ only when the vectors $\xi-\kappa$ and $\kappa$ are located in the opposite side around the origin, that is, $\xi - \kappa \in \mathcal{S}_N^+$ and $\kappa \in \mathcal{S}_N^-$ or $\kappa \in \mathcal{S}_N^+$ and $\xi - \kappa \in \mathcal{S}_N^-$. Then, the angle $\beta$ between two vectors $\xi - \kappa$ and $\kappa$ satisfies 
\[\pi - \tan^{-1}\left(\frac{1}{N-\frac12}\right)<\beta<\pi + \tan^{-1}\left(\frac{1}{N-\frac12}\right).\]
Then, by taking sufficiently large $N$, we make
\begin{equation}\label{eq:k-x.k}
-\frac34 \ge \cos(\beta) = \frac{(\xi -\kappa)\cdot \kappa}{|\xi-\kappa| |\kappa|}.
\end{equation}
Moreover, on each support of $\phi_N(\xi-\kappa)$ and $\phi_N(\kappa)$, we have
\begin{equation}\label{eq:-p}
\begin{aligned}
-p(\xi,\kappa)
=&\frac{(\kappa_1-\xi_1+\kappa_2-\xi_2)}{|\xi-\kappa|}
\frac{(\kappa_1+\kappa_2)}{|\kappa|}
\geq \frac{\left(N-\frac32\right)^2}{1+(\frac12+N)^2}
\geq \frac34,
\end{aligned}
\end{equation}
provided that $N>16$. Thus, for $t:= \frac{1}{100 N }$, by Lemma \ref{lem:cos sin_2D}, \eqref{eq:k-x.k} and \eqref{eq:-p}, we get
\[ \begin{aligned}
        \Vert A_2 (\vec{v}_0) \Vert_{H^{s'}(\R)\times H^{s'}(\R)\times H^{s'}(\R)}
        &\ge \left\Vert \langle \xi\rangle^{s'}\int_0^t  \begin{pmatrix} \widehat{Q}_1 \\ \widehat{Q}_2 \\ \widehat{Q}_3\end{pmatrix} \; ds \right\Vert_{(L^2 \times L^2 \times L^2)(\R^2)}\\
        &\ge \left\Vert \langle \xi\rangle^{s'} \int_0^t  \widehat{Q}_2 ds \right\Vert_{L^2(\R^2)}\\
         &\geq \frac{1}{64}\cdot \frac{3}{4} \left\Vert \langle \xi\rangle^{s'}
				  \dfrac{\xi_1}{(1+d|\xi|^2)}
					\int_0^t  \int_{\R^2}\widehat{\psi}_{N}(\kappa)\widehat{\psi}_{N}(\xi-\kappa)  d\kappa ds \right\Vert_{L^2_{\xi}(S_L)}\\
		&  \gtrsim N^{-2s-1}\left\Vert \langle \xi\rangle^{s'}
		\dfrac{\xi_1}{(1+d|\xi|^2)}
		 \right\Vert_{L^2_{\xi}(S_L)}\\
    \end{aligned}\]
which does not guarantee the uniform boundedness of $\Vert A_2 (\vec{v}_0) \Vert_{H^{s'}(\R^2)\times H^{s'}(\R^2)}$ for $s < -\frac12$. This completes the proof.
\end{proof}

\begin{rem}
Thanks to the symmetric structure of $Q_2$ and $Q_3$ (see Appendix \ref{app:C} for more details), the same conclusion is obtained by taking $Q_3$ as a target instead of $Q_2$ in the proof of  Theorem \ref{thm:2d} \eqref{Thm4}.
	\end{rem}
\subsection{KdV-KdV regime: $b= d = 0$ and $a=c=\frac16$}\footnote{By scaling, we make $a=c=1$ as one dimensional case.}
For $\varrho_K$ is given in \eqref{eq:varrhoK}, from the observation below
\[ | |\xi| \varrho_K(|\xi|) | =|\xi| |1-|\xi|^2|\sim |\xi|^3 ,\]
for $|\xi| \gg 1$, we have
\begin{lem}\label{lem:cos sin_2D_KdV}
	Let $N \gg 1$ be sufficiently large,  $T =\frac{1}{100 N^3}$, and $0 \le s \le t \le T$.  If $|\kappa|, |\xi - \kappa| \sim N$, $|\xi| \sim 1$ and $(\kappa-\xi)\cdot \kappa >0$, then we have
	\begin{equation*}
	\begin{gathered}
	  \frac{\xi\cdot \kappa}{|\xi| |\kappa| }
	\Jbk(t,\xi)  \Jbk(s,\xi-\kappa) \Jak(s,\kappa) 
	+\frac{(\kappa-\xi)\cdot \kappa}{|\xi-\kappa| |\kappa|}
	\frac12 \Jak(t, \xi) \Jak(s,\xi-\kappa)\Jak(s,\kappa) \\
	\geq 
	\frac{1}{16}\left( 	\frac{(\kappa-\xi)\cdot \kappa}{|\xi-\kappa||\kappa|} -\frac{1}{2}\right)
	\end{gathered}
	\end{equation*}
	where $\wh{J_j^K}$, $j=1,2$ are as in \eqref{eq:J2D_K}.
\end{lem}

\begin{proof}
The proof is analogous to the proof of  Lemma \ref{lem:cos sin_2D}, thus we omit the details.
\end{proof}
\subsection*{Proof of Theorem \ref{thm:2d} \eqref{Thm5}}
Suppose that the flow map $\vv_0 \mapsto \vv[\vv_0]$ is continuous in $H^s$, $s < -\frac 32$. Taking the initial data $\vv_0=(\eta_0,u_{01},u_{02})$ as in  \eqref{eq:initial_data}, we obtain
\[\vv_{1}=\begin{pmatrix} \eta_1 \\ u_{11}\\u_{12}  \end{pmatrix} =\mathbf{S}_K(t) \begin{pmatrix} \eta_0 \\ u_{01}\\ u_{02}  \end{pmatrix}=\left(\begin{array}{c}
J^K_1\eta_0+\frac{\partial_{x_1}}{|\nabla|} J^K_2u_{01}+\frac{\partial_{x_2}}{|\nabla|} J^K_2 u_{02} \\
\frac{\partial_{x_1}}{|\nabla|} J^K_2  \eta_0+\frac{-\partial_{x_1}^2}{|\nabla|^2}J^K_1 u_{01} +\frac{-\partial{x_1}\partial_{x_2}}{|\nabla|^2}J^K_1 u_{02}\\
\frac{\partial_{x_2}}{|\nabla|} J^K_2  \eta_0 +\frac{-\partial_{x_1}\partial_{x_2}}{|\nabla|^2}J^K_1 u_{01}+\frac{-\partial_{x_2}^2}{|\nabla|^2}J^K_1 u_{02}
\end{array}\right)
,\]
thus so
\[
A_2(\vec{v_1})=\int_{0}^t \mathbf{S}_K(t-s)
    \left(\begin{array}{c}
                \partial_{x_1}(\eta_{1} u_{11})+\partial_{x_2}(\eta_1 u_{12})    \\
                2^{-1}\partial_{x_1}(u_{11}^2+ u_{12}^2) \\
                2^{-1}\partial_{x_2}(u_{11}^2 +u_{12}^2)
    \end{array}\right)
    =: \int_{0}^t \left(\begin{array}{c}
     Q_1(s) \\
     Q_2(s)\\
     Q_3(s)
    \end{array} \right) ds.
    \]
Analogously to the proof of Theorem \ref{thm:2d} \eqref{Thm4}, we focus on $Q_2$. A direct computation yields that $Q_2$ is decomposed as $Q_2 = Q_{21}+Q_{22}$, where 
    \begin{equation*}
    \begin{aligned}
    \wh{Q}_{21}=&-
    i\xi_1 \int_{\R^2}
    \frac{\xi\cdot \kappa}{|\xi| |\kappa| }
    p(\xi,\kappa) \Jbk(t,\xi)  \Jbk(s,\xi-\kappa) \Jak(s,\kappa) \widehat{\psi}_N(\xi-\kappa) \widehat{\psi}_N(\kappa) d\kappa\\
    \wh{Q}_{22}=& i \xi_1  
    \int_{\R^2} \frac{(\xi-\kappa)\cdot \kappa}{|\xi-\kappa| |\kappa|}
    p(\xi,\kappa)\frac12 \Jak(t, \xi) \Jak(s,\xi-\kappa)\Jak(s,\kappa) \widehat{\psi}_N(\xi-\kappa)  \widehat{\psi}_N (\kappa)d\kappa,
    \end{aligned}
    \end{equation*}
for $p$ given by \eqref{eq:pp}.  Note that all computations in Appendix \ref{app:D} are available for KdV-KdV case (by putting $b=d=0$ and $a=c=1$). Thus, for $t:= \frac{1}{100 N^3 }$, by Lemma \ref{lem:cos sin_2D_KdV} , \eqref{eq:k-x.k} and \eqref{eq:-p}, we get
\[ \begin{aligned}
        \Vert A_2 (\vec{v}_0) \Vert_{H^{s'}(\R)\times H^{s'}(\R)\times H^{s'}(\R)}       &\geq \frac{1}{64}\cdot \frac{3}{4} \left\Vert \langle \xi\rangle^{s'}
				  \xi_1
					\int_0^t  \int_{\R^2}\widehat{\psi}_{N}(\kappa)\widehat{\psi}_{N}(\xi-\kappa)  d\kappa ds \right\Vert_{L^2_{\xi}(S_L)}\\
		&\gtrsim N^{-2s-3}\left\Vert \langle \xi\rangle^{s'} \xi_1
		 \right\Vert_{L^2_{\xi}(S_L)}
    \end{aligned}\]
which does not guarantee the uniform boundedness of $\Vert A_2 (\vec{v}_0) \Vert_{H^{s'}(\R^2)\times H^{s'}(\R^2)}$ for $s < -\frac32$. This completes the proof.

\subsection{BBM-BBM regime: $a=c=0$ and $b=d=1/6$}
We have

\begin{lem}\label{lem:cos sin_2D_BBM}
	Let $N \gg 1$ be sufficiently large, and $T = \frac{1}{1000 }$ and $0 \le s \le t \le T$. If $|\kappa|, |\xi - \kappa| \sim N$ and $|\xi| \sim 1$, then we have
	\begin{equation*}
	\begin{aligned}
	\frac{\xi\cdot \kappa}{|\xi| |\kappa| }
	\Jbb(t,\xi)  \Jbb(s,\xi-\kappa) \Jab(s,\kappa)
	+&
	\frac{(\kappa-\xi)\cdot \kappa}{|\xi-\kappa| |\kappa|}
	\frac12\Jab(t, \xi) \Jab(s,\xi-\kappa)\Jab(s,\kappa)\\
	&\geq \frac{1}{16}\left( 	\frac{(\kappa-\xi)\cdot \kappa}{|\xi-\kappa||\kappa|} -\frac{1}{2}\right),
	\end{aligned}
	\end{equation*}
	where $\varsigma$ and $\wh{J_j^B}$ with $j=1,2$ are as in \eqref{eq:eigenvalues}, \eqref{eq:J2D_BBM} and  respectively.
\end{lem}
\begin{proof}
The proof is analogous to the proof of Lemma \ref{lem:cos sin_2D} with
	\[|\xi| \varrho_B(|\xi|) = \frac{|\xi|}{(1+\frac16|\xi|^2)}
	\sim |\xi|^{-1}
	,\]
for $|\xi| \gg 1$.
\end{proof}

\subsection*{Proof of Theorem \ref{thm:2d} \eqref{Thm3}}
Suppose that the flow map $\vv_0 \mapsto \vv[\vv_0]$ is continuous in $H^s$, $s < 0$. Taking the same initial data $\vv_0=(\eta_0,u_{01},u_{02})$ as in \eqref{eq:initial_data}, we compute
\[\vv_{1}=\begin{pmatrix} \eta_1 \\ u_{11}\\u_{12}  \end{pmatrix} =\mathbf{S}_B(t) \begin{pmatrix} \eta_0 \\ u_{01}\\ u_{02}  \end{pmatrix}
=
\left(\begin{array}{c}
J^B_1\eta_0+\frac{\partial_{x_1}}{|\nabla|} J^B_2u_{01}+\frac{\partial_{x_2}}{|\nabla|} J^B_2 u_{02} \\
\frac{\partial_{x_1}}{|\nabla|} J^B_2  \eta_0+\frac{-\partial_{x_1}^2}{|\nabla|^2}J^B_1 u_{01} +\frac{-\partial{x_1}\partial_{x_2}}{|\nabla|^2}J^B_1 u_{02}\\
\frac{\partial_{x_2}}{|\nabla|} J^B_2  \eta_0 +\frac{-\partial_{x_1}\partial_{x_2}}{|\nabla|^2}J^B_1 u_{01}+\frac{-\partial_{x_2}^2}{|\nabla|^2}J^B_1 u_{02}
\end{array}\right)
,\]
where $J^B_1$ and $J^B_2$ defined by \eqref{eq:J2D_BBM}, thus so
\[
A_2(\vec{v_1})=\int_{0}^t \mathbf{S}_B(t-s)
\left(\begin{array}{c}
(1-\frac16 \Delta)^{-1}[\partial_{x_1}(\eta_{1} u_{11})+\partial_{x_2}(\eta_1 u_{12})]     \\
2^{-1}(1-\frac16\Delta)^{-1}\partial_{x_1}(u_{11}^2+ u_{12}^2) \\
2^{-1}(1-\frac16\Delta)^{-1}\partial_{x_2}(u_{11}^2 +u_{12}^2)
\end{array}\right)
=: \int_{0}^t \left(\begin{array}{c}
Q_1(s) \\
Q_2(s)\\
Q_3(s)
\end{array} \right) ds.
\]
 Note also that all computations in Appendix \ref{app:D} are available for KdV-KdV case (by putting $a=c=0$ and $b=d=1/6$).Thus, $Q_2$ is decomposed as $Q_2 = Q_{21} + Q_{22}$, where 
\[\begin{aligned}
\wh{Q}_{21}=&-
i\xi_1 \int_{\R^2} 
\frac{\xi\cdot \kappa}{|\xi| |\kappa| }
\frac{p(\xi,\kappa)}{(1+\frac16|\xi|^2)} \Jbb(t,\xi)  \Jbb(s,\xi-\kappa) \Jab(s,\kappa) \widehat{\psi}_N(\xi-\kappa) \widehat{\psi}_N(\kappa) d\kappa,\\
\wh{Q}_{22}=& i \xi_1 
\int_{\R^2} \frac{(\xi-\kappa)\cdot \kappa}{|\xi-\kappa| |\kappa|}
\frac{p(\xi,\kappa)}{(1+\frac16 |\xi|^2)}\frac12  \Jab(t, \xi) \Jab(s,\xi-\kappa)\Jab(s,\kappa) \widehat{\psi}_N(\xi-\kappa)  \widehat{\psi}_N (\kappa)d\kappa,
\end{aligned}\]
for $p$ given by \eqref{eq:pp}.
Thus, for $t:= \frac{1}{1000}$, by Lemma \ref{lem:cos sin_2D_BBM}, \eqref{eq:k-x.k} and \eqref{eq:-p}, we get
\[ \begin{aligned}
\Vert A_2 (\vec{v}_0) \Vert_{H^{s'}(\R)\times H^{s'}(\R)\times H^{s'}(\R)} &\geq \frac{1}{64}\cdot \frac{3}{4} \left\Vert \langle \xi\rangle^{s'}
				  \xi_1
					\int_0^t  \int_{\R^2}\widehat{\psi}_{N}(\kappa)\widehat{\psi}_{N}(\xi-\kappa)  d\kappa ds \right\Vert_{L^2_{\xi}(S_L)}\\
&\gtrsim N^{-2s} \left\Vert \langle \xi\rangle^{s'} \xi_1
\right\Vert_{L^2_{\xi}(S_L)},
\end{aligned}\]
which does not guarantee the uniform boundedness of $\Vert A_2 (\vec{v}_0) \Vert_{H^{s'}(\R^2)\times H^{s'}(\R^2)}$ for $s < 0$,  this ends the proof.
\begin{flushright}
	\qedsymbol
\end{flushright}

\appendix

\section{Local Well-Posedness}\label{app:A}
This section briefly shows the local well-posedness of \eqref{eq:abcd} and \eqref{eq:2D_ABCD} including BBM-BBM case, but not KdV-KdV case. This result may not be optimal except for BBM-BBM case. The well-known well-posedness theorem is given by
\begin{thm}\label{thm:LWP}
Let $n=1,2$. Fix $s\geq 0$. For any $(u_0,v_0)\in H^s(\R^n)\times H^s(\R^n)$ , there exists a $T(u_0,v_0) > 0$ and a  unique solution $(u,v)\in X^{s}_{T}$ (for a suitable solution space $X_T^s$) of the initial value problem \eqref{eq:abcd}. The maximal existence time $T=T_s$ for the solution has the property that
	\[
	  T_s\geq \dfrac{ C_s}{\|(u_0,v_0) \|_{H^s(\R)\times H^s(\R)}}
	\]
 where the positive constant $C_s$ depends only on $s$.
\end{thm}

It is well-known that Theorem \ref{thm:LWP} immediately follows from multilinear estimates, thus in what follows, we only focus on bilinear estimate (see Section \ref{sec:LWP} below). 

\subsection{Notations}\label{sec:notations}
We define Bessel and Riesz potentials ($J^s$ and $D^s$, respectively) of order $-s$, $s \in \R$, as Fourier multipliers by
\[J^s f (x) := \mathcal F^{-1}\left((1+|\xi|^2)^{\frac{s}{2}} \widehat{f}\right) \quad \mbox{and} \quad D^sf(x) := \mathcal F^{-1} \left( |\xi|^s \widehat{f}\right).\]
In particular, $\sqrt{-\Delta} = D^1 = D$ is the Fourier multiplier of the symbol $|\xi|$.

\subsection{Littlewood-Paley Decomposition}\label{LPD}
This section devotes to explaining the Littlewood-Paley decomposition, which is an useful way to improve the bilinear estimate for the local well-posedness theory. As well-known, the Littlewood-Paley decomposition is a particular way to write a single function as a superposition of a countably infinite family of functions of varying frequencies.

\medskip

Let $\varphi(\xi)$ be a real-valued radially symmetric bump function on $\R^n$ with the support $\{\xi \in \R^n :  |\xi| \le 2\}$ which is identical to $1$ on the set $\{\xi \in \R^n : |\xi| \le 1\}$ and is decreasing on $\{\xi \in \R^n : 1 \le |\xi| \le 2\}$. Define a dyadic number $N \in 2^{\Z_{\ge 0}}$ of the form $N = 2^k$, $k \in \Z_{\ge 0}$. Let denote $\varphi_1 = \varphi$ and define
\[\varphi_N(\xi) = \varphi\left(\frac{\xi}{N}\right) - \varphi\left(\frac{2\xi}{N}\right), \quad N \ge 2.\]
By construction, the sequence of $\varphi_N$ satisfies
\[\sum_{N \ge 1 \; : \;  dyadic} \varphi_N (\xi) \equiv 1.\]
We simply write $\sum_{N \ge 1}$ by dropping "dyadic". This provides a typical \emph{partition of unity} which allows to define the projection operator (one of the so-called \emph{Littlewood-Paley} projection operator) on $L^2(\R^n)$ by $P_Nf(x) = \mathcal F^{-1} \left(\varphi_N(\xi) \widehat{f} \right) (x)$. Using the projection operators, one decompose any function $f$ in $L^2(\R)$ as
\[f = \sum_{N \ge 1} P_N f.\]
We sometimes denote $P_Nf$ by simply $f_N$. Note that $f_N$ belongs to any Sobolev space $H^s(\R^n)$, $s \ge 0$ (or smooth), whenever $f \in L^2(\R^n)$. The following lemma is well-known \emph{Bernstein inequality}, which is to upgrade low Lebesgue integrability to high Lebesgue integrability with the price of some powers of $N$:
\begin{lem}[Bernstein's inequalities]
Let $f \in L^2(\R^n)$, $1 \le p \le q \le \infty$, and $s \ge 0$. Then, we have
\begin{equation}\label{eq:Bernstein_1}
\norm{D^{\pm s}f_N}_{L^p} \sim N^{\pm s} \norm{f_N}_{L^p}
\end{equation}
and
\begin{equation}\label{eq:Bernstein_2}
\norm{f_N}_{L^q} \lesssim N^{\frac{n}{p} - \frac{n}{q}} \norm{f_N}_{L^p}.
\end{equation}
The implicit constants in both \eqref{eq:Bernstein_1} and \eqref{eq:Bernstein_2} depend only on $s, n, p$ and $q$.
\end{lem}

\subsection{Bilinear estimates}\label{sec:LWP}

\begin{lem}\label{lem:bilinear_0}
Let $f, g \in L^2(\R^n)$. Then, we have
\[\norm{P_1(fg)}_{L^2} \lesssim \norm{f}_{L^2}\norm{g}_{L^2}.\]
\end{lem}

\begin{proof}
The proof follows from $|(\wh{f}\ast\wh{g})(\xi)| \lesssim \norm{f}_{L^2}\norm{g}_{L^2}$ and $\int_{|\xi| \le 1} \; d\xi \lesssim 1$. 
\end{proof}

\begin{lem}[Refined bilinear estimate]\label{lem:bilinear}
Let $s \ge \frac{n-2}{2}$ and $f, g \in H^s(\R^n)$. Then, we have
\begin{equation}\label{eq:bilinear}
\norm{J^{-1}D(fg)}_{H^s} \lesssim \norm{f}_{H^s}\norm{g}_{H^s}.
\end{equation}
\end{lem}

\begin{proof}
From the duality argument, it suffices for the left-hand side of \eqref{eq:bilinear} to estimate
\[\int_{\R^2} \left(J^{s-1}D(fg) \right) w \; dx,\]
where $w \in L^2$ with $\norm{w}_{L^2} =1$. We make the Littlewood-Paley decomposition of $f,g$ and $w$ as 
\[f = \sum_{N_1 \ge 1} f_{N_1}, \quad g = \sum_{N_2 \ge 1}g_{N_2} \quad \mbox{and} \quad w = \sum_{N \ge 1} w_N,\]
respectively. Without loss of generality, we may assume $N_1 \le N_2$. By Lemma \ref{lem:bilinear_0}, we are now reduced to establishing
\begin{equation}\label{eq:bi2}
\sum_{N > 1}\sum_{N_1, N_2 \ge 1}N^{-1+s}\int f_{N_1}g_{N_2} w_N \; dx \lesssim \norm{f}_{H^s}\norm{g}_{H^s}.
\end{equation}
We separate \eqref{eq:bi2} into two cases :  $N \sim N_2 \gtrsim N_1$ and $N_2 \sim N_1 \gg N$.

\medskip

\textbf{(Case I.)} $N \sim N_2 \gtrsim N_1$. Using H\"older's and Bernstein's \eqref{eq:Bernstein_2} inequalities, one has
\[\int f_{N_1}g_{N_2} w_N \; dx \le \norm{f_{N_1}}_{L^{\infty}}\norm{g_{N_2}}_{L^2}\norm{w_N}_{L^2} \lesssim N_1^{\frac n2}\norm{f_{N_1}}_{L^2}\norm{g_{N_2}}_{L^2}\norm{w_N}_{L^2}.\]
With this, we further reduce \eqref{eq:bi2} to
\begin{equation}\label{eq:bi3}
\sum_{N > 1}\sum_{N_2 \sim N}\sum_{N_1 \lesssim N_2} N^{-1+s}N_2^{-s}N_1^{\frac n2-s} \norm{f_{N_1}}_{H^s}\norm{g_{N_2}}_{H^s}\norm{w_N}_{L^2} \lesssim \norm{f}_{H^s}\norm{g}_{H^s},
\end{equation}
thanks to \eqref{eq:Bernstein_1}. We denote the multiplier in the left-hand side of \eqref{eq:bi3} by $m(N,N_1,N_2)$, namely, $m(N,N_1,N_2) = N^{-1+s}N_2^{-s}N_1^{1-s}$. Note that the number of $N_2$ is finite, and then
\[\sum_{N}\sum_{N_2 \sim N} \norm{g_{N_2}}_{L^2}^2 \sim \sum_{N}\norm{g_N}_{L^2}^2.\]
With this observation, Cauchy-Schwarz inequality yields
\[\mbox{LHS of } \eqref{eq:bi3} \lesssim \left(\sup_{N > 1} \sum_{\substack{1 \le N_1 \lesssim N_2 \\
N_2 \sim N}} m^2(N,N_1,N_2) \right)^{\frac12}\norm{f}_{H^s}\norm{g}_{H^s}\norm{w}_{L^2}.\]
Hence, it remains to prove 
\begin{equation}\label{eq:bi4}
\sup_{N > 1} \sum_{N_2 \sim N}\sum_{1 \le N_1 \lesssim N_2} N^{-2+2s}N_2^{-2s}N_1^{n-2s} \lesssim 1.
\end{equation}
This is sometimes referred as \emph{Schur's test}, see, for instance, \cite[Lemma 3.11]{Tao2001} for more details. Then, one gets if $s > \frac n2$
\[\sum_{N_2 \sim N}\sum_{1 \le N_1 \lesssim N_2} N^{-2+2s}N_2^{-2s}N_1^{n-2s} \lesssim \sum_{N_2 \sim N}N^{-2+2s}N_2^{-2s} \lesssim N^{-2},\]
otherwise ($s \le \frac n2$),
\[\sum_{N_2 \sim N}\sum_{1 \le N_1 \lesssim N_2} N^{-2+2s}N_2^{-2s}N_1^{n-2s} \lesssim \sum_{N_2 \sim N}N^{-2+2s}N_2^{n-4s} \lesssim N^{n-2-2s}.\]
Thus, \eqref{eq:bi4} holds true if $s \ge \frac{n-2}{2}$.

\medskip

\textbf{(Case II.)} $N_2 \sim N_1 \gg N$. Analogously, using H\"older's and Bernstein's \eqref{eq:Bernstein_2} inequalities, one obtains
\[\int f_{N_1}g_{N_2} w_N \; dx \le \norm{f_{N_1}}_{L^2}\norm{g_{N_2}}_{L^2}\norm{w_N}_{L^{\infty}} \lesssim N^{\frac n2}\norm{f_{N_1}}_{L^2}\norm{g_{N_2}}_{L^2}\norm{w_N}_{L^2},\]
which ensures \eqref{eq:bilinear} provided that
\begin{equation}\label{eq:bi5}
\sup_{N_2 > 1} \sum_{N_1 \sim N_2}\sum_{1 < N \lesssim N_1} N_2^{-2s}N_1^{-2s}N^{n-2+2s} \lesssim 1.
\end{equation}
One immediately obtains \eqref{eq:bi5} for $s \ge \frac{n-2}{2}$, we thus complete the proof.
\end{proof}

\begin{rem}
Lemma \ref{lem:bilinear} slightly improves the bilinear estimates by Grisvard \cite{Grisvard}, particularly, validity of the bilinear estimates in $H^{\frac{n-2}{2}}$. This seems to facilitate the global well-posedness in the energy space for $4$-dimensional problem.   
\end{rem}

\begin{rem}\label{rem:b,d constant}
Replacing the Bessel potential $J^{-1}$ by $(1-b\Delta)^{-1}$ or $(1-d\Delta)^{-1}$ in Lemma \ref{lem:bilinear} affects only the constant in the left-hand side, precisely, the implicit constant should depend on $b$ or $d$.
\end{rem}

\begin{rem}\label{rem:LWP}
The  standard Picard iteration method immediately assures the (smooth) local well-posedness in $H^{s}(\R^n)$, $s = \max(0, \frac{n-2}{2})$. Thus, BBM-BBM case is optimal in the sense that the flow map is analytic.
\end{rem}

\section{Decomposition of $(Q_1,Q_2)$ on the generic regime: one-dimensional case}\label{app:1D_decomposition}
This section provides a precise computation of $Q_j$, $j=1,2$, for one-dimensional case. Recall \eqref{eq:A_2}
\[
A_2(\vec{v_1})=\int_{0}^t S(t-s)
\left(\begin{array}{c}
(1-b\partial_x^{2})^{-1}\partial_{x}(\eta_{1} u_{1})     \\
2^{-1}(1-d\partial_x^2)^{-1}\partial_{x}(u_{1}^2) \\
\end{array}\right)
=: \int_{0}^t \left(\begin{array}{c}
Q_1 \\
Q_2\\
\end{array} \right) ds,
\]
where the linear propagator $S$ is given in \eqref{eq:semigroup_l}. A direct computation yields
\begin{equation}\label{eq:Q1-2}
\begin{aligned}
\widehat{Q}_1
    =&~{} i \xi \left[ \frac{1}{(1+b\xi^2)}\La(t-s,\xi)\widehat{\eta_1 u_1}(s,\xi) - \frac{h(\xi)}{2(1+d\xi^2)}\Lb(t-s,\xi)\widehat{u_1^2}(s, \xi) \right]\\
\widehat{Q}_2
    =&~{} i \xi \left[-\frac{1}{h(\xi)(1+b\xi^2)}\Lb(t-s,\xi)\widehat{\eta_1 u_1}(s,\xi) + \frac{1}{2(1+d\xi^2)}\La(t-s,\xi)\widehat{u_1^2}(s, \xi) \right]
    \end{aligned}
\end{equation}
Inserting the initial data $\vec{v}_{0}=(\eta_0, u_0)=(0,\phi)$ into \eqref{eq:v1}, we have
\[
\mathcal{F}\left(\begin{array}{c}
\eta_1\\
u_{1}
\end{array} \right)=\mathcal{F}\left(S(t)\left(\begin{array}{c}
\eta_0\\
u_{0}
\end{array}\right) \right)
=
\begin{pmatrix}
-h(\xi)\widehat{L}_2(t,\xi) \widehat{\phi} \\
 \widehat{L}_1(t,\xi) \widehat{\phi}
\end{pmatrix} .
\]
With this, a direct computation gives
\begin{equation*}
\begin{aligned}
\widehat{\eta_1 u_{1}}(\xi)
= \int_{\R}\widehat{\eta_{1}}(\xi_1)\widehat{u_{1}}(\xi-\xi_1)d\xi_1
=-\int_{\R}
h(\xi_1)\widehat{L}_2(t,\xi_1) \widehat{L}_1(t,\xi-\xi_1) \widehat{\phi}(\xi_1) \widehat{\phi}(\xi-\xi_1)
\end{aligned}
\end{equation*}
and
\begin{equation*}
\begin{aligned}
\widehat{u_{1} u_{1}}(\xi)
= \int_{\R}\widehat{u_{1}}(\xi-\xi_1) \widehat{u_{1}}(\xi_1)d\xi_1
&=\int_{\R}
 \widehat{L}_1(t,\xi_1) \widehat{L}_1(t,\xi-\xi_1)
 \widehat{\phi}(\xi_1)  \widehat{\phi} (\xi-\xi_1)d\xi_1.
\end{aligned}
\end{equation*}
Thus, by inserting them into \eqref{eq:Q1-2}, we conclude that
\begin{equation*}
\begin{aligned}
\widehat{Q}_1
    = -i \xi \Bigg[&~{}\frac{1}{(1+b\xi^2)} \int_{\R} h(\xi_1)\La(t-s,\xi) \Lb(s,\xi_1)
    \La(s,\xi-\xi_1) \widehat{\phi}(\xi_1)  \widehat{\phi}(\xi-\xi_1) d\xi_1 \\
    & +\frac{h(\xi)}{2(1+d\xi^2)} \int_{\R}  \Lb(t-s,\xi)\La(s,\xi_1) 
    \La(s,\xi-\xi_1)\widehat{\phi}(\xi_1) \widehat{\phi}(\xi-\xi_1)d\xi_1 \Bigg]
    \end{aligned}
\end{equation*}

and
\begin{equation*}
\begin{aligned}
\widehat{Q}_2 =
    i \xi \Bigg[&~{}\frac{1}{h(\xi)(1+b\xi^2)} \int_{\R} h(\xi_1)\Lb(t-s,\xi) \Lb(s,\xi_1)
    \La(s,\xi-\xi_1) \widehat{\phi}(\xi_1)  \widehat{\phi}(\xi-\xi_1) d\xi_1 \\
    & +\frac{1}{2(1+d\xi^2)} \int_{\R}  \La(t-s,\xi)\La(s,\xi_1) 
    \La(s,\xi-\xi_1)\widehat{\phi}(\xi_1) \widehat{\phi}(\xi-\xi_1)d\xi_1 \Bigg].
    \end{aligned}
\end{equation*}

\section{Decomposition of $(Q_1,Q_2,Q_3)$: generic case} \label{app:C}
This section provides a precise computation of $Q_j$, $j=1,2,3$, for two-dimensional case. Recall \eqref{eq:semigroup_2d}
   \begin{equation}\label{eq:semigroup_2d_S}
\begin{aligned}
& \mathcal F  \left(\mathbf{S}(t) \begin{pmatrix} f \\ g \\ h\end{pmatrix}\right)
 =&~{}  \begin{pmatrix} 
 \Ja(t, \xi) \widehat{f}  +
    \varsigma(|\xi|)\frac{i}{|\xi|}\Jb(t,\xi) \big[\xi_1\widehat{g} +  \xi_2\wh{h}\big]
  \\
 \frac{i\xi_1}{\varsigma(|\xi|)|\xi|}\Jb(t,\xi) \widehat{f}  + 
 \frac{\xi_1}{|\xi|^2} \Ja(t, \xi) \big[ \xi_1 \widehat{g}+ \xi_2  \wh{h} \big]\\
 \frac{i\xi_2}{\varsigma(|\xi|)|\xi|}\Jb(t,\xi) \widehat{f}  + 
 \frac{\xi_2}{|\xi|^2} \Ja(t, \xi) \big[ \xi_1 \widehat{g}+ \xi_2  \wh{h} \big]
  \end{pmatrix} 
 ,
\end{aligned}
\end{equation}
where $J_j$, $j=1,2$, are given in \eqref{eq:J2D}. Inserting the initial data $\vec{v}_{0}=(\eta_0, u_{01}, u_{02})=(0, \phi, \phi) = (f,g,h)$ into \eqref{eq:semigroup_2d_S}, we have
\begin{equation}\label{eq:eta1u11u12}
\mathcal{F} \begin{pmatrix} 
\eta_1\\
u_{11}\\
u_{12}
  \end{pmatrix} 
=
 \begin{pmatrix} 
\varsigma(|\xi|)\frac{i(\xi_1 + \xi_2)}{|\xi|}\Jb(t,\xi)\wh{\phi}\\
  \frac{\xi_1(\xi_1 + \xi_2)}{|\xi|^2} \Ja(t, \xi) \wh{\phi}\\
 \frac{\xi_2(\xi_1 + \xi_2)}{|\xi|^2} \Ja(t, \xi)  \wh{\phi}
  \end{pmatrix} .
\end{equation}
Recall the integral part in \eqref{eq:Duhamel_2d}
\[
A_2(\vec{v_1})=\int_{0}^t \mathbf{S}(t-s)
\left(\begin{array}{c}
(1-b\Delta)^{-1}[\partial_{x_1}(\eta_{1} u_{11})+\partial_{x_2}(\eta_1 u_{12})]     \\
2^{-1}(1-d\Delta)^{-1}\partial_{x_1}(u_{11}^2+ u_{12}^2) \\
2^{-1}(1-d\Delta)^{-1}\partial_{x_2}(u_{11}^2 +u_{12}^2)
\end{array}\right)
=: \int_{0}^t \left(\begin{array}{c}
Q_1 \\
Q_2\\
Q_3
\end{array} \right) ds.
\]
A direct computation on $Q_j$, $j=1,2,3$, yields
\begin{equation}\label{eq:Q1-3}
\begin{aligned}
\wh{Q}_1=&  \frac{i}{(1+b|\xi|^2)} \Ja(t,\xi)\left[\xi_1(\wh{\eta_{1} u_{11}})+\xi_2(\wh{\eta_1 u_{12}})\right]  
-
 \frac{\varsigma(|\xi|) |\xi|}{2(1+d|\xi|^2)} \Jb(t,\xi) \left(\wh{u_{11}^2} +\wh{u_{12}^2}\right) ,\\
\wh{Q}_2=&-
\frac{\xi_1}{\varsigma(|\xi|)|\xi|}\frac{\Jb(t,\xi)}{(1+b|\xi|^2)}  \left[\xi_1(\wh{\eta_{1} u_{11}})+\xi_2(\wh{\eta_1 u_{12}})\right]  
+ 
  \frac{i \xi_1}{2(1+d|\xi|^2)} \Ja(t, \xi)\left(\wh{u_{11}^2} +\wh{u_{12}^2}\right) ,\\
\wh{Q}_3=&-\frac{\xi_2}{\varsigma(|\xi|)|\xi|} \frac{\Jb(t,\xi)}{(1+b|\xi|^2)} 
 \left[\xi_1(\wh{\eta_{1} u_{11}})+\xi_2(\wh{\eta_1 u_{12}})\right]  + 
 \frac{i\xi_2}{2(1+d|\xi|^2)} \Ja(t, \xi) \left(\wh{u_{11}^2} +\wh{u_{12}^2}\right).
\end{aligned}
\end{equation}
With \eqref{eq:eta1u11u12}, we first compute nonlinear interactions
\[\begin{aligned}
\widehat{\eta_1 u_{11}}(\xi)
= i\int_{\R^2}
\varsigma(|\kappa|)p(\xi,k)\frac{(\xi_1 - \kappa_1)}{|\xi - \kappa|} \Jb(t,\kappa)\Ja(t, \xi - \kappa) \wh{\phi}(\kappa)\wh{\phi}(\xi-\kappa) \; d\kappa,
\end{aligned}\]
\[\begin{aligned}
\widehat{\eta_1 u_{12}}(\xi)
= i\int_{\R^2}
\varsigma(|\kappa|)p(\xi,k)\frac{(\xi_2 - \kappa_2)}{|\xi - \kappa|} \Jb(t,\kappa)\Ja(t, \xi - \kappa) \wh{\phi}(\kappa)\wh{\phi}(\xi-\kappa) \; d\kappa,
\end{aligned}\]
\[\begin{aligned}
\widehat{u_{11}^2}(\xi)
= \int_{\R^2}
p(\xi,k)\frac{\kappa_1}{|\kappa|}\frac{(\xi_1 - \kappa_1)}{|\xi - \kappa|} \Ja(t,\kappa)\Ja(t, \xi - \kappa) \wh{\phi}(\kappa)\wh{\phi}(\xi-\kappa) \; d\kappa,
\end{aligned}\]
and
\[\begin{aligned}
\widehat{u_{12}^2}(\xi)
= \int_{\R^2}
p(\xi,k)\frac{\kappa_2}{|\kappa|}\frac{(\xi_2 - \kappa_2)}{|\xi - \kappa|} \Ja(t,\kappa)\Ja(t, \xi - \kappa) \wh{\phi}(\kappa)\wh{\phi}(\xi-\kappa) \; d\kappa,
\end{aligned}\]
where
	\[	p(\xi,\kappa)= \frac{(\xi_1+\xi_2-\kappa_1-\kappa_2)}{|\xi-\kappa|}
		\frac{(\kappa_1+\kappa_2)}{|\kappa|}.\]
Inserting them into \eqref{eq:Q1-3}, we obtain

\[\begin{aligned}
\wh{Q}_1=&  
-\frac{1}{(1+b|\xi|^2)}
 \int_{\R^2}\varsigma(|\kappa|)p(\xi,\kappa) \frac{\xi \cdot (\xi-\kappa)}{|\xi - \kappa|}\Ja(t, \xi) \Jb(s,\xi-\kappa) \Ja(s,\kappa) \widehat{\phi}(\kappa) \widehat{\phi}(\xi-\kappa) d\kappa\\
&-
\frac{\varsigma(|\xi|) |\xi|}{2(1+d|\xi|^2)} 
\int_{\R^2}p(\xi,\kappa) 
\frac{\kappa \cdot (\xi-\kappa)}{|\xi-\kappa||\kappa| }
	\Jb(t,\xi)  \Ja(s,\xi-\kappa)\Ja(s,\kappa) \widehat{\phi} (\kappa) \widehat{\phi}(\xi-\kappa) d\kappa,
\end{aligned}\]
\[\begin{aligned}
\wh{Q}_2=&-
\frac{i\xi_1}{(1+b|\xi|^2)} 
 \int_{\R^2} p(\xi,\kappa)\frac{\varsigma(|\kappa|)}{\varsigma(|\xi|)} \frac{\xi \cdot (\xi-\kappa)}{|\xi||\xi - \kappa|}\Jb(t,\xi) \Jb(s,\xi-\kappa) \Ja(s,\kappa) \widehat{\phi}(\kappa) \widehat{\phi}(\xi-\kappa) d\kappa\\
&+ 
  \frac{i \xi_1}{2(1+d|\xi|^2)} 
\int_{\R^2} \frac{(\xi-\kappa)\cdot \kappa}{|\xi-\kappa| |\kappa|}
p(\xi,\kappa) \Ja(t, \xi) \Ja(s,\xi-\kappa)\Ja(s,\kappa)  \widehat{\phi}(\kappa) \widehat{\phi}(\xi-\kappa) d\kappa
\end{aligned}\]
and
\[\begin{aligned}
\wh{Q}_3=&-
\frac{i\xi_2}{(1+b|\xi|^2)} 
 \int_{\R^2} p(\xi,\kappa)\frac{\varsigma(|\kappa|)}{\varsigma(|\xi|)} \frac{\xi \cdot (\xi-\kappa)}{|\xi||\xi - \kappa|}\Jb(t,\xi) \Jb(s,\xi-\kappa) \Ja(s,\kappa) \widehat{\phi}(\kappa) \widehat{\phi}(\xi-\kappa) d\kappa\\
&+ 
  \frac{i \xi_2}{2(1+d|\xi|^2)} 
\int_{\R^2} \frac{(\xi-\kappa)\cdot \kappa}{|\xi-\kappa| |\kappa|}
p(\xi,\kappa) \Ja(t, \xi) \Ja(s,\xi-\kappa)\Ja(s,\kappa)  \widehat{\phi}(\kappa) \widehat{\phi}(\xi-\kappa) d\kappa.
\end{aligned}\]

\section{Decomposition of $(Q_1,Q_2,Q_3)$: $a=c$ and $b=d \geq 0$} \label{app:D}
This section provides a precise computation of $Q_j$, $j=1,2,3$, for two-dimensional case when $a=c$ and $b=d\geq0$. Note that the decomposition here is valid for both KdV-KdV and BBM-BBM cases. Recall \eqref{eq:semigroup_2d_ab}
\[\begin{aligned}
& \mathcal F  \left(\mathbf{S}_{ab}(t) \begin{pmatrix} f \\ g \\ h\end{pmatrix}\right)
=&~{}  \begin{pmatrix} 
\Ja^{ab}(t, \xi) \widehat{f}  +
\frac{i}{|\xi|}\Jb^{ab}(t,\xi) \big[\xi_1\widehat{g} +  \xi_2\wh{h}\big]
\\
\frac{i\xi_1}{|\xi|}\Jb^{ab}(t,\xi) \widehat{f}  + 
\frac{\xi_1}{|\xi|^2} \Ja^{ab}(t, \xi) \big[ \xi_1 \widehat{g}+ \xi_2  \wh{h} \big]\\
\frac{i\xi_2}{|\xi|}\Jb^{ab}(t,\xi) \widehat{f}  + 
\frac{\xi_2}{|\xi|^2} \Ja^{ab}(t, \xi) \big[ \xi_1 \widehat{g}+ \xi_2  \wh{h} \big]
\end{pmatrix} 
,
\end{aligned}\]
where $J_j^{ab}$, $j=1,2$, are given in \eqref{eq:J2D_ab}. Inserting the initial data $\vec{v}_{0}=(\eta_0, u_{01}, u_{02})=(0, \phi, \phi) = (f,g,h)$ into $\vec{v}_{1}=\mathbf{S}_{ab}\vec{v}_{0}$ 
, we have
\begin{equation}\label{eq:eta1u11u12_ab}
\mathcal{F} \begin{pmatrix} 
\eta_1\\
u_{11}\\
u_{12}
\end{pmatrix} 
=
\begin{pmatrix} 
\frac{i(\xi_1 + \xi_2)}{|\xi|}\Jb^{ab}(t,\xi)\wh{\phi}\\
\frac{\xi_1(\xi_1 + \xi_2)}{|\xi|^2} \Ja^{ab}(t, \xi) \wh{\phi}\\
\frac{\xi_2(\xi_1 + \xi_2)}{|\xi|^2} \Ja^{ab}(t, \xi)  \wh{\phi}
\end{pmatrix} .
\end{equation}
Recall the integral part in \eqref{eq:Duhamel_2d_ab}
\[
A_2(\vec{v_1})=\int_{0}^t \mathbf{S}_{ab}(t-s)
\left(\begin{array}{c}
(1-b\Delta)^{-1}[\partial_{x_1}(\eta_{1} u_{11})+\partial_{x_2}(\eta_1 u_{12})]     \\
2^{-1}(1-b\Delta)^{-1}\partial_{x_1}(u_{11}^2+ u_{12}^2) \\
2^{-1}(1-b\Delta)^{-1}\partial_{x_2}(u_{11}^2 +u_{12}^2)
\end{array}\right)
=: \int_{0}^t \left(\begin{array}{c}
Q_1 \\
Q_2\\
Q_3
\end{array} \right) ds.
\]
A direct computation on $Q_j$, $j=1,2,3$, yields
\begin{equation}\label{eq:Q1-3_ab}
\begin{aligned}
\wh{Q}_1=&  \frac{i}{(1+b|\xi|^2)} \Ja^{ab}(t,\xi) \left[\xi_1(\wh{\eta_{1} u_{11}})+\xi_2(\wh{\eta_1 u_{12}})\right]  
-
\frac{|\xi|}{2(1+b|\xi|^2)} \Jb^{ab}(t,\xi) \left(\wh{u_{11}^2} +\wh{u_{12}^2}\right) ,\\
\wh{Q}_2=&-
\frac{\xi_1}{|\xi|}\frac{\Jb^{ab}(t,\xi)}{(1+b|\xi|^2)}  \left[\xi_1(\wh{\eta_{1} u_{11}})+\xi_2(\wh{\eta_1 u_{12}})\right]  
+ 
\frac{i \xi_1}{2(1+b|\xi|^2)} \Ja^{ab}(t, \xi)\left(\wh{u_{11}^2} +\wh{u_{12}^2}\right) ,\\
\wh{Q}_3=&-\frac{\xi_2}{|\xi|} \frac{\Jb^{ab}(t,\xi)}{(1+b|\xi|^2)} 
\left[\xi_1(\wh{\eta_{1} u_{11}})+\xi_2(\wh{\eta_1 u_{12}})\right]  + 
\frac{i\xi_2}{2(1+b|\xi|^2)} \Ja^{ab}(t, \xi) \left(\wh{u_{11}^2} +\wh{u_{12}^2}\right).
\end{aligned}
\end{equation}
With \eqref{eq:eta1u11u12_ab}, we first compute nonlinear interactions
\[\begin{aligned}
\widehat{\eta_1 u_{11}}(\xi)
= i\int_{\R^2}
p(\xi,k)\frac{(\xi_1 - \kappa_1)}{|\xi - \kappa|} \Jb^{ab}(t,\kappa)\Ja^{ab}(t, \xi - \kappa) \wh{\phi}(\kappa)\wh{\phi}(\xi-\kappa) \; d\kappa,
\end{aligned}\]
\[\begin{aligned}
\widehat{\eta_1 u_{12}}(\xi)
= i\int_{\R^2}
p(\xi,k)\frac{(\xi_2 - \kappa_2)}{|\xi - \kappa|} \Jb^{ab}(t,\kappa)\Ja^{ab}(t, \xi - \kappa) \wh{\phi}(\kappa)\wh{\phi}(\xi-\kappa) \; d\kappa,
\end{aligned}\]
\[\begin{aligned}
\widehat{u_{11}^2}(\xi)
= \int_{\R^2}
p(\xi,k)\frac{\kappa_1}{|\kappa|}\frac{(\xi_1 - \kappa_1)}{|\xi - \kappa|} \Ja^{ab}(t,\kappa)\Ja^{ab}(t, \xi - \kappa) \wh{\phi}(\kappa)\wh{\phi}(\xi-\kappa) \; d\kappa,
\end{aligned}\]
and
\[\begin{aligned}
\widehat{u_{12}^2}(\xi)
= \int_{\R^2}
p(\xi,k)\frac{\kappa_2}{|\kappa|}\frac{(\xi_2 - \kappa_2)}{|\xi - \kappa|} \Ja^{ab}(t,\kappa)\Ja^{ab}(t, \xi - \kappa) \wh{\phi}(\kappa)\wh{\phi}(\xi-\kappa) \; d\kappa,
\end{aligned}\]
where
\[	p(\xi,\kappa)= \frac{(\xi_1+\xi_2-\kappa_1-\kappa_2)}{|\xi-\kappa|}
\frac{(\kappa_1+\kappa_2)}{|\kappa|}.\]
Inserting them into \eqref{eq:Q1-3_ab}, we obtain

\[\begin{aligned}
\wh{Q}_1=&  
-\frac{1}{(1+b|\xi|^2)}
\int_{\R^2}p(\xi,\kappa) \frac{\xi \cdot (\xi-\kappa)}{|\xi - \kappa|}\Ja^{ab}(t, \xi) \Jb^{ab}(s,\xi-\kappa) \Ja^{ab}(s,\kappa) \widehat{\phi}(\kappa) \widehat{\phi}(\xi-\kappa) d\kappa\\
&-
\frac{ |\xi|}{2(1+b|\xi|^2)} 
\int_{\R^2}p(\xi,\kappa) 
\frac{\kappa \cdot (\xi-\kappa)}{|\xi-\kappa||\kappa| }
\Jb^{ab}(t,\xi)  \Ja^{ab}(s,\xi-\kappa)\Ja^{ab}(s,\kappa) \widehat{\phi} (\kappa) \widehat{\phi}(\xi-\kappa) d\kappa,
\end{aligned}\]
\[\begin{aligned}
\wh{Q}_2=&-
\frac{i\xi_1}{(1+b|\xi|^2)} 
\int_{\R^2} p(\xi,\kappa) \frac{\xi \cdot (\xi-\kappa)}{|\xi||\xi - \kappa|}\Jb^{ab}(t,\xi) \Jb^{ab}(s,\xi-\kappa) \Ja^{ab}(s,\kappa) \widehat{\phi}(\kappa) \widehat{\phi}(\xi-\kappa) d\kappa\\
&+ 
\frac{i \xi_1}{2(1+b|\xi|^2)} 
\int_{\R^2} \frac{(\xi-\kappa)\cdot \kappa}{|\xi-\kappa| |\kappa|}
p(\xi,\kappa) \Ja^{ab}(t, \xi) \Ja^{ab}(s,\xi-\kappa)\Ja^{ab}(s,\kappa)  \widehat{\phi}(\kappa) \widehat{\phi}(\xi-\kappa) d\kappa
\end{aligned}\]
and
\[\begin{aligned}
\wh{Q}_3=&-
\frac{i\xi_2}{(1+b|\xi|^2)} 
\int_{\R^2} p(\xi,\kappa) \frac{\xi \cdot (\xi-\kappa)}{|\xi||\xi - \kappa|}\Jb^{ab}(t,\xi) \Jb^{ab}(s,\xi-\kappa) \Ja^{ab}(s,\kappa) \widehat{\phi}(\kappa) \widehat{\phi}(\xi-\kappa) d\kappa\\
&+ 
\frac{i \xi_2}{2(1+b|\xi|^2)} 
\int_{\R^2}  p(\xi,\kappa) \frac{(\xi-\kappa)\cdot \kappa}{|\xi-\kappa| |\kappa|}
 \Ja^{ab}(t, \xi) \Ja^{ab}(s,\xi-\kappa)\Ja^{ab}(s,\kappa)  \widehat{\phi}(\kappa) \widehat{\phi}(\xi-\kappa) d\kappa.
\end{aligned}\]

\end{document}